\documentclass[11pt,reqno, a4paper]{amsart}
\usepackage[T1]{fontenc}
\usepackage[utf8]{inputenc}
\usepackage{lmodern}
\usepackage{amsmath,amssymb,amsthm,mathtools}
\usepackage[margin=1in]{geometry}
\usepackage{microtype}
\usepackage{cite}
\usepackage[colorlinks=true,linkcolor=blue,citecolor=blue,urlcolor=blue]{hyperref}

\newtheorem{theorem}{Theorem}[section]
\newtheorem{proposition}[theorem]{Proposition}
\newtheorem{lemma}[theorem]{Lemma}
\newtheorem{corollary}[theorem]{Corollary}
\theoremstyle{definition}

\newtheorem{example}[theorem]{Example}
\theoremstyle{remark}
\newtheorem{remark}[theorem]{Remark}

\newcommand{\R}{\mathbb R}
\newcommand{\C}{\mathbb C}
\newcommand{\F}{\mathbb F}
\newcommand{\one}{\mathbf 1}
\newcommand{\id}{\mathrm{id}}
\newcommand{\Sp}{S_0^+}
\newcommand{\qn}{q_{n,p}}
\newcommand{\Hp}{H_{p}}
\newcommand{\rp}{\rho_p}
\newcommand{\norm}[1]{\| #1\|}
\newcommand{\supnorm}[1]{\norm{#1}_{\infty}}
\newcommand{\abs}[1]{\left\lvert #1\right\rvert}
\DeclareMathOperator{\supp}{supp}
\DeclareMathOperator{\spanop}{span}

\title[Sphere isometries of smooth functions]{Extensions from positive and full unit spheres of $C^\infty([0,1])$ with compatible Fréchet metrics}
\author{Min-Ruei Lin}
\address{Department of Applied Mathematics, National Sun Yat-sen University, Kaohsiung 80424, Taiwan; and Graduate School of Science and Technology, Niigata University, Niigata 950--2181, Japan}
\email{m082030021@student.nsysu.edu.tw}

\begin{document}

\begin{abstract}
We study how the order of positivity determines the form of a linear isometric extension from the positive unit sphere of $C^\infty([0,1],\mathbb F)$, where $\mathbb F=\mathbb R$ or $\mathbb C$.
At order zero, positivity means that $f(t)\in[0,\infty)$ for every $t\in[0,1]$.
At order $m\ge1$, we require $f^{(m)}(t)\in[0,\infty)$ for every $t\in[0,1]$ and $f^{(j)}(x_0)\in[0,\infty)$ for $0\le j<m$, where $x_0$ is a fixed anchor.
For the compatible Fr\'echet metrics $\rho_{m,p}$ studied here, we prove
that every surjective isometry of the positive unit sphere extends uniquely to a surjective $\mathbb F$-linear isometry of the whole space.
This holds for all $1\le p\le\infty$ at order zero and for $p=\infty$ under an explicit anchor condition at higher orders.
At order zero, the extension either fixes every function or acts by reflection of the interval. Reflection is possible exactly when $p=\infty$ or $x_0=1/2$.
At higher orders, the extension acts on the $m$th derivative by the identity or reflection and may permute the lower-order derivative values at the anchor.
A counterexample shows that, for a different compatible metric, reflection can preserve distances on the positive unit sphere without admitting a real-linear isometric extension.

We also classify the surjective isometries of the full unit sphere at order zero and prove that each extends uniquely to a surjective real-linear isometry of the whole space.
\end{abstract}

\maketitle

\section{Introduction}\label{sec:intro}
Tingley's problem~\cite{Tingley} asks whether every surjective isometry between the unit spheres of real Banach spaces extends to a surjective real-linear isometry of the ambient spaces; see Peralta~\cite{PeraltaSurvey} for a survey in the setting of operator algebras.
Recently, OpenAI~\cite{OpenAITingley} announced an affirmative solution for arbitrary real Banach spaces in a preprint.
Alongside the classical problem, extension questions for smaller subsets and for spaces beyond the Banach setting have developed in their own right.
The positive unit sphere is a natural smaller domain for such a question.
If the positive cone generates the real ambient space, then the positive unit sphere spans that space and hence determines any possible real-linear extension uniquely.
The issue is whether an isometry of this generating subset necessarily admits an isometric extension.
Leung, Ng, and Wong~\cite{LeungNgWong} obtained positive-sphere extension results for $L^p$ spaces with $1\le p\le\infty$ and continuous-function spaces, while Ezumi, Lin, and Miura~\cite{EzumiLinMiura} established corresponding results for continuous functions vanishing at infinity.
For spaces involving differentiation, positivity can also incorporate derivative information.
Related choices of order in absolutely continuous and Sobolev function spaces are studied in~\cite{LinAC,LinSobolev}.

The present paper develops this question on $C^\infty([0,1],\mathbb F)$, with $\mathbb F\in\{\mathbb R,\mathbb C\}$, through a family of positivity conditions indexed by derivative order.
At order zero, positivity is defined pointwise.
At each positive order, the derivative of that order is pointwise nonnegative, and the lower-order derivatives are nonnegative at a fixed anchor.
The higher-order setting thus replaces pointwise positivity of the function by positivity of its $m$th derivative, together with nonnegativity of the lower-order derivatives at the anchor.
These conditions are naturally expressed through the anchored-coordinate representation.
Our aim is to determine how this passage from functions to derivative data is reflected in the existence, uniqueness, and form of isometric extensions.
The usual smooth topology introduces an essential difficulty: it is not induced by a norm.
We therefore equip the space with compatible translation-invariant metrics that are unbounded and have convex balls.
These metrics retain useful features of norm geometry while accounting for derivatives of every order.
However, different compatible metrics can give different sphere geometries, so the metric choice is part of the extension problem.
This places our work within the broader study of sphere extensions beyond Banach spaces~\cite{Tan,Fu} and of metrics and linear isometries on Fr\'echet spaces~\cite{Muller,ChalendarOgerPartington}.

Our positive-sphere theorem describes all extensions in terms of the lower-order derivatives at the anchor and the $m$th derivative.
At order zero, only the function component remains.
The proof recovers these coordinates from distances on the positive sphere, then uses the derivative terms in the metric to determine the admissible interval maps.
The full-sphere theorem concerns the order-zero metric over both scalar fields.
It gives uniqueness among real-linear extensions and yields a classification of the ambient isometries.

The metric assumptions alone do not guarantee these conclusions.
We construct another compatible metric with the stated geometric properties for which a positive-sphere isometry has no real-linear isometric extension (see Subsection~\ref{subsec:fail_extension}).
Thus, throughout the family, the extension problem depends on how the chosen metric records the function and its derivative data.
The classification results and the counterexample together make this dependence explicit.

The paper is organized as follows.
Section~\ref{sec:notation} introduces the notation and states the principal results.
Section~\ref{sec:metric} establishes the metric properties used below.
Section~\ref{sec:positive} proves the positive-sphere theorem for the family indexed by $m\ge0$ and gives a counterexample for a different metric.
Section~\ref{sec:full} treats the full order-zero sphere and classifies the global isometries.
The appendices contain the metric ball formulas and the anchored-coordinate representation.

\section{Anchored coordinates, metrics, and main results}\label{sec:notation}
Consider $C^\infty([0,1],\F)$ with $\F\in\{\R,\C\}$.
A function is smooth on $[0,1]$ if it is smooth on $(0,1)$ and its derivatives of every order extend continuously to the endpoints.
We equip $ C^\infty([0,1], \F)$ with the usual Fr\'echet topology generated by
\[
 P_n(f):=\max_{0\le j\le n}\supnorm{f^{(j)}},\qquad n\ge0.
\]
Thus convergence means uniform convergence of every derivative.
We write
\[
 C^\infty([0,1], \C)_\R= C^\infty([0,1], \R)\oplus i C^\infty([0,1], \R),
\]
where $ C^\infty([0,1], \C)_\R$ denotes the $\C$-valued space with scalars restricted to $\R$.
For $\F\in\{\R,\C\}$, the space $C^\infty([0,1],\F)$ is regarded as a vector space over $\F$, equipped with its usual Fr\'echet topology.
Unless otherwise specified, a linear map on $ C^\infty([0,1], \F)$ is $\F$-linear.
The notation $ C^\infty([0,1], \F)$ always refers to this specific smooth-function space.
We denote by $\one$ the constant function with value one.

Fix an \textit{anchor} $x_0\in[0,1]$ and $1\le p\le\infty$. For $m\ge0$ write
\[
 J_m f=((f^{(j)}(x_0))_{0\le j<m}, f^{(m)})\qquad(f\in C^\infty([0,1],\F)).
\]
When $m=0$, the finite block is empty and we identify $J_0f$ with $f$.
Define
\begin{equation}\label{eq:anchored}
 \qn(f):=
 \begin{cases}
 \displaystyle\left(\sum_{j=0}^{n-1}|f^{(j)}(x_0)|^p+
                   \supnorm{f^{(n)}}^p\right)^{1/p},\quad&p<\infty,\\[6pt]
 \displaystyle\max\left\{\max_{0\le j<n}|f^{(j)}(x_0)|,\,
                   \supnorm{f^{(n)}}\right\},\quad&p=\infty.
 \end{cases}
\end{equation}
Put $q_0=q_{0,p}=\supnorm{\cdot}$,
\begin{equation}\label{eq:metric}
 H_p(h):=\sup_{n\ge1}2^{-n}\frac{q_{n,p}(h)}{1+q_{n,p}(h)},
 \qquad \rho_{m,p}(f,g):=\max\{q_{m,p}(f-g),H_p(f-g)\}.
\end{equation}
We suppress the anchor in the notation. Let $\mathcal T$ denote the \textit{usual Fr\'echet topology}, $\mathcal T_\infty$ the uniform topology, and $\mathcal T_{\rho_p}$ the topology induced by $\rho_p=\rho_{0,p}$.
We also refer to $\mathcal T$ as the \textit{smooth topology}.

In $C^\infty([0,1],\F)$, identify $[0,\infty)$ with the nonnegative real scalars in $\F$. Define
\[
 K_m:=\left\{f\in C^\infty([0,1],\F):
 \begin{array}{ll}
 f^{(j)}(x_0)\in[0,\infty) & (0\le j<m),\\
 f^{(m)}(t)\in[0,\infty) & (t\in[0,1])
 \end{array}\right\},
\]
\[
 \qquad S_m^+:=\{f\in K_m:\rho_{m,p}(f,0)=1\}.
\]
These conditions give the same cone under the natural inclusion of the $\R$-valued space into the $\C$-valued space.
In particular, for $\F=\C$ and $f\in K_m$, the function $v=\operatorname{Im}f$ satisfies $v^{(m)}=0$ and $v^{(j)}(x_0)=0$ for $0\le j<m$. Taylor's formula gives $v=0$ for $m\ge1$, and for $m=0$ this follows from pointwise positivity.
At order zero, $K_0$ is the pointwise nonnegative cone.

Let $r(x)=1-x$. We call the identity admissible for every $p,x_0$, and call $r$ admissible if $p=\infty$ or $x_0=1/2$.
Define
\begin{equation}\label{eq:L_and_B}
L_{x_0}:=\max\{x_0,1-x_0\}\qquad \text{and}\qquad B_j(L_{x_0}):=\sum_{k=0}^j \frac{L_{x_0}^k}{k!}\quad (j\ge0).
\end{equation}
We state the main theorems of this paper.

\begin{theorem}\label{thm:main}
Fix an integer $m\ge0$ and suppose that one of the following holds.
\begin{enumerate}
\item $m=0$ and $1\le p\le\infty$.
\item $m\ge1$, $p=\infty$, and $B_{m-1}(L_{x_0})\le2$.
\end{enumerate}
Let $\F\in\{\R,\C\}$ and $\Delta^+:S_m^+\to S_m^+$ be a surjective $\rho_{m,p}$-isometry.
Then $\Delta^+$ admits a unique $\F$-linear surjective $\rho_{m,p}$-isometric extension $\widetilde{\Delta^+}:C^\infty([0,1],\F)\to C^\infty([0,1],\F)$.
More precisely, there exist a unique permutation $\pi$ of $\{0,\ldots, m-1\}$ and a unique admissible $\sigma\in\{\id, r\}$ such that
\begin{equation}\label{eq:finite-midpoint-extension}
 J_m(\widetilde{\Delta^+}f)=\bigl((f^{(\pi(j))}(x_0))_{0\le j<m},\ f^{(m)}\circ\sigma\bigr).
\end{equation}
At $m=0$, the formula~\eqref{eq:finite-midpoint-extension} reads $\widetilde{\Delta^+}f=f\circ\sigma$.
Conversely, every such pair $(\pi,\sigma)$ gives a surjective ambient isometry whose restriction is a surjective $\rho_{m,p}$-isometry of $S_m^+$.
\end{theorem}
We write the full unit sphere of $C^\infty([0,1], \F)$ as
\[
S_\F:=\{f\in C^\infty([0,1],\F): \rho_{0,p}(f,0) =1\}
\]
\begin{theorem}\label{thm:full}
Let $\F\in\{\R,\C\}$ and let $\Delta:S_{\F}\to S_{\F}$ be a surjective $\rp$-isometry.
Then $\Delta$ admits a unique $\R$-linear surjective $\rp$-isometric extension $\widetilde{\Delta}:
 C^\infty([0,1],\F)\to C^\infty([0,1],\F)$.
More precisely, there exist a constant $\lambda\in\mathbb{T}_{\F}:=\{z\in \F:|z|=1\}$, an admissible map $\sigma:[0,1]\to[0,1]$, and a map $\kappa:\F\to\F$ such that
\begin{equation}\label{eq:fullrepresentation}
 \widetilde{\Delta}f
 =\lambda\,\kappa(f\circ\sigma)
 \qquad(f\in C^\infty([0,1],\F)).
\end{equation}
Here $\kappa$ is the identity when $\F=\R$ and either the identity or complex conjugation when $\F=\C$, acting pointwise on functions.
In the complex case, the extension is complex linear or conjugate linear, according to the choice of $\kappa$.
Conversely, for every choice of $\lambda$, $\sigma$, and $\kappa$ as above, the restriction of $f\mapsto\lambda\,\kappa(f\circ\sigma)$ to $S_{\F}$ is a surjective $\rp$-isometry.
\end{theorem}
Sections~\ref{sec:positive} and~\ref{sec:full} prove Theorems~\ref{thm:main} and~\ref{thm:full}, respectively.

The following spanning identities give uniqueness of a linear extension whenever it exists. For $m\ge1$, Proposition~\ref{prop:anchored-coordinates} identifies $K_m$ with $\R_+^m\oplus C^\infty([0,1],\R)_+$ in anchored coordinates.
Thus $K_m$ is closed, convex, and pointed.
When the anchor is the left endpoint, repeated integration yields $K_m\subseteq K_0$. For $x_0>0$, the function $t\mapsto t-x_0$ lies in $K_m$ and is negative at zero.

\begin{proposition}\label{prop:finite-jet-span}
Let $m\ge0$ and $1\le p\le\infty$.
Then
\[
 K_m-K_m=C^\infty([0,1],\R),\qquad
 \spanop_{\R}S_m^+=C^\infty([0,1],\R),\qquad
 \spanop_{\C}S_m^+=C^\infty([0,1],\C).
\]
In particular, $\spanop_{\R}S_m^+$ is a proper real
subspace of $C^\infty([0,1],\C)_{\R}$.
\end{proposition}
\begin{proof}
Let $f\in C^\infty([0,1],\R)$. For $m=0$, choose
$M\ge\supnorm f$. Then $f+M\one$ and $M\one$ are pointwise
nonnegative, so $f=(f+M\one)-M\one\in K_0-K_0$.

For $m\ge1$, write $J_mf=(b_0,\ldots,b_{m-1},g)$ and choose $M\ge\max\{|b_0|,\ldots,|b_{m-1}|,\supnorm g\}$.
Both tuples on the right-hand side of
\[
 J_mf=
 (b_0+M,\ldots,b_{m-1}+M,g+M\one)
 -(M,\ldots,M,M\one)
\]
have nonnegative real entries in the finite block and a pointwise nonnegative function in the final coordinate.
By Proposition~\ref{prop:anchored-coordinates}, applying $J_m^{-1}$ to these tuples gives two members of $K_m$ whose difference is $f$.
Thus $C^\infty([0,1],\R)\subseteq K_m-K_m$ for every $m\ge0$.
The reverse is immediate, and hence $K_m-K_m=C^\infty([0,1],\R)$.

For each nonzero $h\in K_m$,
\[
 h=q_{m,p}(h)\frac{h}{q_{m,p}(h)},
 \qquad
 \frac{h}{q_{m,p}(h)}\in S_m^+.
\]
It follows that $K_m$ is contained in $\spanop_{\R}S_m^+$.
Since $K_m-K_m=C^\infty([0,1],\R)$ and $S_m^+$ consists of $\R$-valued functions, we obtain $\spanop_{\R}S_m^+=C^\infty([0,1],\R)$.

Now let $f\in C^\infty([0,1],\C)$.
Write $f=\operatorname{Re}f+i\operatorname{Im}f$.
Both $\operatorname{Re}f$ and $\operatorname{Im}f$ belong to $\spanop_{\R}S_m^+$ by the real spanning identity.
Hence $\spanop_{\C}S_m^+=C^\infty([0,1],\C)$.
Finally, $i\one$ belongs to $C^\infty([0,1],\C)_{\R}$ and lies outside the $\R$-valued space. Thus $\spanop_{\R}S_m^+$ is a proper real subspace of $C^\infty([0,1],\C)_{\R}$.
\end{proof}

\section{Metric and convex geometry}\label{sec:metric}
Fix $m\ge0$ and $1\le p\le\infty$. We show that $\rho_{m,p}$ induces the usual Fr\'echet topology, describe its balls and its behavior on rays, and recover the base norm $q_{m,p}$ from metric distances.

\begin{proposition}\label{prop:norms}
For every $n\ge1$, $\qn$ is a norm.
Moreover, $(q_0,q_{1,p},q_{2,p},\ldots)$ generates $\mathcal{T}$.
\end{proposition}
\begin{proof}
Fix $n\geq 1$. We first show that $\qn$ is a norm.
Absolute homogeneity is immediate from the definition.
For $1\leq p<\infty$, the inequalities
\[
 |(f+g)^{(j)}(x_0)|
 \le |f^{(j)}(x_0)|+|g^{(j)}(x_0)|,
 \qquad
 \supnorm{(f+g)^{(n)}}
 \le \supnorm{f^{(n)}}+\supnorm{g^{(n)}}
\]
and Minkowski's inequality in $\ell_p^{n+1}$ give $\qn(f+g)\le \qn(f)+\qn(g)$.
For $p=\infty$, the triangle inequality follows by taking the maximum of the corresponding estimates.
Suppose that $\qn(f)=0$.
Then
\[
f^{(j)}(x_0)=0 \quad (0\leq j<n),
\qquad
f^{(n)}=0.
\]
Since $f^{(n)}=0$, the function $f$ is a polynomial of degree at most $n-1$. 
All its Taylor coefficients at $x_0$ vanish, and hence $f=0$. 
Thus $\qn$ is a norm.

We next prove the norm equivalence.
By the definition of $P_n$, we have
\[
|f^{(j)}(x_0)|\leq P_n(f)
\quad (0\leq j<n),
\qquad
\supnorm{f^{(n)}}\leq P_n(f).
\]
Consequently, for $1\leq p<\infty$,
\begin{align*}
\qn(f)^p=\sum_{j=0}^{n-1}|f^{(j)}(x_0)|^p
  +\supnorm{f^{(n)}}^p\leq (n+1)P_n(f)^p.
\end{align*}
For $p=\infty$, the same inequalities give $q_{n,\infty}(f)\leq P_n(f)$.
Therefore,
\[
\frac{1}{C_{n,p}}\qn(f)\leq P_n(f),
\]
where $C_{n,p}=(n+1)^{1/p}$ for $1\le p<\infty$, and $C_{n,\infty}=1$.

We now establish the reverse estimate.
Fix $0\leq k<n$.
Taylor's formula with integral remainder gives, for every $x\in[0,1]$,
\begin{align*}
f^{(k)}(x)=\sum_{j=k}^{n-1}
  \frac{f^{(j)}(x_0)}{(j-k)!}
  (x-x_0)^{j-k} + \frac{1}{(n-k-1)!}\int_{x_0}^{x}(x-t)^{n-k-1}f^{(n)}(t)\,dt.
\end{align*}
Taking absolute values and estimating the integral remainder, we obtain
\begin{align*}
|f^{(k)}(x)|\leq
\sum_{j=k}^{n-1}
\frac{|f^{(j)}(x_0)|}{(j-k)!}
|x-x_0|^{j-k} + \frac{\supnorm{f^{(n)}}}{(n-k)!}
|x-x_0|^{n-k}.
\end{align*}
Since $x,x_0\in[0,1]$, we have $|x-x_0|\leq 1$.
Moreover, $|f^{(j)}(x_0)|\leq \qn(f)$ and $\supnorm{f^{(n)}}\leq \qn(f)$.
It follows that
\[
|f^{(k)}(x)|\leq \qn(f)
\left(
\sum_{j=k}^{n-1}\frac{1}{(j-k)!}
+\frac{1}{(n-k)!}
\right)=\qn(f)
\sum_{\ell=0}^{n-k}\frac{1}{\ell!}
\leq e\,\qn(f),
\]
where $e$ is Euler's number.
Since the preceding inequality holds for every $x\in[0,1]$, we have
\begin{equation}\label{eq:bound1}
\supnorm{f^{(k)}}\leq e\,\qn(f)
\qquad (0\leq k<n).
\end{equation}
For $k=n$, the definition of $\qn$ gives
\begin{equation}\label{eq:bound2}
\supnorm{f^{(n)}}\leq \qn(f).
\end{equation}
Combining inequalities~\eqref{eq:bound1} and \eqref{eq:bound2} with $P_n$, we conclude that $P_n(f)\leq e\,\qn(f)$.

Thus, for every $n\geq 1$,
\[
\frac{1}{C_{n,p}}\qn(f)
\leq P_n(f)
\leq e\,\qn(f).
\]
That is, $\qn$ and $P_n$ are equivalent norms for each fixed $n$. 
In particular, since $q_0=P_0$, convergence with respect to the family $(q_0,q_{1,p},q_{2,p},\ldots)$ is equivalent to convergence with respect to $(P_n)_{n\geq 0}$.
Hence $(q_0,q_{1,p},q_{2,p},\ldots)$ generates the usual Fr\'echet topology $\mathcal{T}$.
\end{proof}

We next characterize the metric properties of $\rho_{m,p}$.
Put $\phi(t)=t/(1+t)$ for $t\ge0$.
Then $\Hp(h) = \sup_{n\ge1}2^{-n}\phi(\qn(h))$ by~\eqref{eq:metric}.
First, we establish the following lemma.

\begin{lemma}\label{lem:tail}
For every $h\in C^\infty([0,1],\F)$ and integer $N\ge0$,
\[
0\le\sup_{n>N}2^{-n}\phi(\qn(h))
\le2^{-(N+1)}.
\]
If $h\ne0$, the supremum defining $\Hp(h)$ is attained at a finite index.
In particular, $0\le\Hp(h)<1/2$ for every $h$.
\end{lemma}
\begin{proof}
Fix $h\in C^\infty([0,1],\F)$, and put $a_n=2^{-n}\phi(\qn(h))$ for $n\ge1$.
Since $0\le\phi(t)<1$ for every $t\ge0$, we have $0\le a_n<2^{-n}$.
Thus, for every $N\ge0$,
\[
0\le\sup_{n>N}a_n
\le\sup_{n>N}2^{-n}
=2^{-(N+1)}.
\]

Suppose $h\ne0$.
Since $q_{1,p}$ is a norm, we have $a_1>0$.
Choose $N\ge1$ such that $2^{-(N+1)}<a_1$.
Then
\[
\sup_{n>N}a_n
\le2^{-(N+1)}<a_1
\le\max_{1\le n\le N}a_n.
\]
Consequently, we have $\Hp(h)=\max_{1\le n\le N}a_n$, so the supremum is attained at a finite index.
Since $a_n<2^{-n}\le1/2$ for every $n\ge1$, it follows that $0\le\Hp(h)<1/2$.
The case $h=0$ is immediate.
\end{proof}

Recall that $C^\infty([0,1],\F)$ is complete in its usual Fr\'echet topology $\mathcal T$.
Consequently, every translation-invariant metric $d$ inducing $\mathcal T$ is complete.
To see this, let $(f_k)$ be a $d$-Cauchy sequence.
For every $n\ge0$ and $\varepsilon>0$, compatibility gives $\delta>0$ such that
\[
 d(h,0)<\delta\quad\Longrightarrow\quad P_n(h)<\varepsilon.
\]
By translation invariance, $d(f_k-f_\ell,0)=d(f_k,f_\ell)$.
Thus $(f_k)$ is Cauchy with respect to every $P_n$.
Completeness of the usual Fr\'echet space yields $f\in C^\infty([0,1],\F)$ such that $f_k\to f$ in $\mathcal T$.
Since $d$ induces $\mathcal T$, we also have $d(f_k,f)\to0$, proving completeness.

\begin{theorem}\label{thm:metric}
For every $m\ge0$ and $1\le p\le\infty$, $\rho_{m,p}$ is a complete, unbounded, translation-invariant metric inducing the usual Fr\'echet topology.
\end{theorem}
\begin{proof}
For $a,b\ge0$ we have
\[
 \phi(a)+\phi(b)-\phi(a+b)
 =\frac{ab(a+b+2)}{(1+a)(1+b)(1+a+b)}\ge0 .
\]
Since $\phi$ is increasing, each $2^{-n}\phi(q_{n,p}(f-g))$ satisfies the triangle inequality.
The supremum over $n$ and the maximum with $q_{m,p}(f-g)$ preserve this inequality.
The base norm gives positive definiteness and unboundedness. Symmetry and translation invariance follow from the definition.

If $\rho_{m,p}(f_k,f)\to0$, then for each $n\ge1$ we have $\phi(q_{n,p}(f_k-f))\le2^n\rho_{m,p}(f_k,f)\to0$.
Proposition~\ref{prop:norms} therefore gives convergence in the smooth topology.
Conversely, smooth convergence implies convergence in $q_{m,p}$ and in each $q_{n,p}$.
For fixed $N$, the maximum of the first $N$ weighted terms tends to zero, and the remaining supremum is at most $2^{-(N+1)}$ uniformly in $k$ by Lemma~\ref{lem:tail}.
This proves compatibility.
Completeness follows from the preceding observation about translation-invariant metrics.
\end{proof}

\begin{proposition}\label{prop:rays}
For $f\ne0$, the function $t\mapsto\rho_{m,p}(tf,0)$ is continuous and strictly increasing from zero to infinity on $[0,\infty)$.
\end{proposition}
\begin{proof}
Fix $f\ne0$ and put
\[
 R_f(t):=\rho_{m,p}(tf,0)
 =\max\left\{tq_{m,p}(f),
       \sup_{n\ge1}2^{-n}\phi\bigl(tq_{n,p}(f)\bigr)\right\},
 \qquad t\ge0.
\]
By Theorem~\ref{thm:metric}, $\rho_{m,p}$ induces the smooth topology $\mathcal T$.
Since scalar multiplication is continuous in this topology, $R_f$ is continuous on $[0,\infty)$.

To prove strict monotonicity, let $0<s<t$.
By Lemma~\ref{lem:tail}, the supremum defining $H_p(sf)$ is attained at a finite index.
Hence $R_f(s)$ equals either $sq_{m,p}(f)$ or $2^{-n}\phi\bigl(sq_{n,p}(f)\bigr)$ for some $n\ge1$.
The term that attains $R_f(s)$ is strictly larger at $t$, because $q_{m,p}(f)>0$, $q_{n,p}(f)>0$ for every $n\ge1$, and $\phi$ is strictly increasing.
Therefore $R_f(t)>R_f(s)$.
Moreover, $R_f(0)=0$ and $R_f(t)\ge tq_{m,p}(f)>0$ for every $t>0$.
Thus $R_f$ is strictly increasing on $[0,\infty)$.

Finally,
\[
 R_f(t)\ge tq_{m,p}(f)\longrightarrow\infty
 \qquad (t\to\infty),
\]
since $q_{m,p}(f)>0$.
\end{proof}

\begin{proposition}\label{prop:large}
Every open or closed $\rho_{m,p}$-ball is convex.
Moreover, for all $f,g\in C^\infty([0,1],\F)$,
\[
 \rho_{m,p}(f,g)\ge\frac12
 \quad\Longleftrightarrow\quad
 q_{m,p}(f-g)\ge\frac12,
\]
and whenever these conditions hold, $\rho_{m,p}(f,g)=q_{m,p}(f-g)$.
\end{proposition}
\begin{proof}
Let $f,g\in C^\infty([0,1],\F)$ and $0\le t\le1$.
Since $q_{n,p}$ is a norm for every $n\ge1$,
\[
 q_{n,p}(tf+(1-t)g)
 \le t q_{n,p}(f)+(1-t)q_{n,p}(g)
 \le\max\{q_{n,p}(f),q_{n,p}(g)\}.
\]
The same inequality holds for the base norm $q_{m,p}$, including $q_{0,p}=q_0$.
Since $\phi$ is increasing, it follows that
\[
 H_p(tf+(1-t)g)\le\max\{H_p(f),H_p(g)\}.
\]
Consequently,
\[
 \rho_{m,p}(tf+(1-t)g,0)
 \le\max\{\rho_{m,p}(f,0),\rho_{m,p}(g,0)\}.
\]
Thus every open or closed $\rho_{m,p}$-ball centered at zero is convex.
Translation invariance gives convexity for arbitrary centers.

By Lemma~\ref{lem:tail}, $H_p(f-g)<1/2$.
Hence $\rho_{m,p}(f,g)=\max\{q_{m,p}(f-g),H_p(f-g)\}$ is at least $1/2$ exactly when $q_{m,p}(f-g)\ge1/2$, in which case it equals $q_{m,p}(f-g)$.
\end{proof}
Taking $g=0$ in Proposition~\ref{prop:large}, we obtain
\[
 S_m^+=\{f\in K_m:q_{m,p}(f)=1\}.
\]
In particular, $S_0^+$ consists of the nonnegative smooth functions of supremum norm one.

The ball formulas and Minkowski functionals are given in Appendix~\ref{app:balls}.
The large-distance identity also determines the base norm intrinsically.
\begin{theorem}\label{thm:recovery}
For $f,g\in C^\infty([0,1],\F)$,
\begin{equation}\label{eq:base-recovery}
 q_{m,p}(f-g)=
 \sup_{\substack{z:\ \rho_{m,p}(f,z)\ge1/2\\ \rho_{m,p}(g,z)\ge1/2}}
 \bigl|\rho_{m,p}(f,z)-\rho_{m,p}(g,z)\bigr|.
\end{equation}
\end{theorem}
\begin{proof}
Fix $f,g\in C^\infty([0,1],\F)$.
Let $z\in C^\infty([0,1],\F)$ satisfy $\rho_{m,p}(f,z)\ge1/2$ and $\rho_{m,p}(g,z)\ge1/2$.
By Proposition~\ref{prop:large}, we have $\rho_{m,p}(f,z)=q_{m,p}(f-z)$ and $\rho_{m,p}(g,z)=q_{m,p}(g-z)$.
The reverse triangle inequality for $q_{m,p}$ gives
\begin{align*}
 \bigl|\rho_{m,p}(f,z)-\rho_{m,p}(g,z)\bigr|
 &=\bigl|q_{m,p}(f-z)-q_{m,p}(g-z)\bigr|\\
 &\le q_{m,p}\bigl((f-z)-(g-z)\bigr)
 =q_{m,p}(f-g).
\end{align*}
Taking the supremum over all such $z$ gives the upper bound in~\eqref{eq:base-recovery}.

Suppose first that $f\ne g$ and put $h=f-g$.
Since $q_{m,p}(h)>0$, choose $t>0$ such that $tq_{m,p}(h)\ge1/2$. Set $z_0=f+th$.
We then obtain $f-z_0=-th$ and $g-z_0=-(t+1)h$, and hence
\[
 q_{m,p}(f-z_0)=tq_{m,p}(h)\ge\frac12,
 \qquad
 q_{m,p}(g-z_0)=(t+1)q_{m,p}(h)\ge\frac12.
\]
By Proposition~\ref{prop:large},
\[
 \rho_{m,p}(f,z_0)=tq_{m,p}(h),
 \qquad
 \rho_{m,p}(g,z_0)=(t+1)q_{m,p}(h).
\]
Thus $z_0$ satisfies both conditions in the supremum of~\eqref{eq:base-recovery}.
Moreover,
\begin{align*}
 \bigl|\rho_{m,p}(f,z_0)-\rho_{m,p}(g,z_0)\bigr|
 =(t+1)q_{m,p}(h)-tq_{m,p}(h)
 =q_{m,p}(h)=q_{m,p}(f-g).
\end{align*}
This proves equality in~\eqref{eq:base-recovery} when $f\ne g$.

If $f=g$, take $z=f+\one$.
Since $q_{m,p}(\one)=1$, Proposition~\ref{prop:large} gives $\rho_{m,p}(f,z)=\rho_{m,p}(g,z)=1$.
The indexing set is therefore nonempty, and every term in the supremum is zero.
Hence~\eqref{eq:base-recovery} also holds when $f=g$.
\end{proof}

\begin{corollary}\label{cor:affine}
Every surjective global $\rho_{m,p}$-isometry $T:C^\infty([0,1],\F)\to C^\infty([0,1],\F)$ preserves $q_{m,p}$-distances and is affine over $\R$.
\end{corollary}
\begin{proof}
Let $T$ be a surjective $\rho_{m,p}$-isometry and fix $f,g\in C^\infty([0,1],\F)$.
Since $T$ is surjective and preserves $\rho_{m,p}$-distances, the substitution $w=Tz$ carries the set of $z$ satisfying $\rho_{m,p}(f,z)\ge1/2$ and $\rho_{m,p}(g,z)\ge1/2$ bijectively onto the set of $w$ satisfying $\rho_{m,p}(Tf,w)\ge1/2$ and $\rho_{m,p}(Tg,w)\ge1/2$.
Thus Theorem~\ref{thm:recovery} gives
\begin{align*}
 q_{m,p}(f-g)
 &=\sup_{\substack{z\in C^\infty([0,1],\F)\\
        \rho_{m,p}(f,z)\ge1/2,\ \rho_{m,p}(g,z)\ge1/2}}
   \bigl|\rho_{m,p}(f,z)-\rho_{m,p}(g,z)\bigr|\\
 &=\sup_{\substack{w\in C^\infty([0,1],\F)\\
        \rho_{m,p}(Tf,w)\ge1/2,\ \rho_{m,p}(Tg,w)\ge1/2}}
   \bigl|\rho_{m,p}(Tf,w)-\rho_{m,p}(Tg,w)\bigr|
 =q_{m,p}(Tf-Tg).
\end{align*}
Therefore $T$ is a surjective isometry for the norm $q_{m,p}$.

Regard $C^\infty([0,1],\F)$ as a real normed space equipped with $q_{m,p}$, regardless of whether $\F=\R$ or $\F=\C$.
The Mazur--Ulam theorem applies without a completeness assumption~\cite{Vaisala}.
Hence $T$ is affine over $\R$.
\end{proof}
At $m=0$, formula~\eqref{eq:base-recovery} recovers the supremum distance from $\rho_p$.
The same argument applies to isometries between two order-zero metrics having different anchors or exponents.

\section{Isometries of positive spheres}\label{sec:positive}
We prove Theorem~\ref{thm:main} by recovering the anchored coordinates from distances on $S_m^+$ and using them to determine the linear extension of surjective isometries.
We then give a counterexample showing how extension depends on the metric. For a different compatible metric, reflection preserves distances on $S_0^+$ but admits no $\R$-linear isometric extension to the ambient space.

\subsection{Classification and linear extension of positive-sphere isometries}\label{subsec:positive-classification}
Unless otherwise specified, we assume either $m=0$ and $1\le p\le\infty$, or $m\ge1$ and $p=\infty$.
The anchor is arbitrary.
Let
\[
 Y_m:=\{a_0,\ldots,a_{m-1}\}\sqcup[0,1]
\]
be the disjoint union of the interval and $m$ isolated points, with $Y_0=[0,1]$.
View the coordinate space from Proposition~\ref{prop:anchored-coordinates} as
\[
 \mathcal A_m:=\R^m\oplus C^\infty([0,1],\R)\subseteq C(Y_m,\R),
\]
where smoothness is required on the interval component.
The jet cone and base norm take the form
\[
 J_m(K_m)=\mathcal A_m^+,
 \qquad q_{m,p}(f)=\|J_mf\|_{C(Y_m)}.
\]
Put
\[
 \Sigma_m^+:=\{F\in\mathcal A_m:F\ge0,
                    \ \|F\|_{C(Y_m)}=1\},
 \qquad
 d_m(F,G):=\rho_{m,p}(J_m^{-1}F,J_m^{-1}G).
\]
By Proposition~\ref{prop:large}, for every \(f\in K_m\) we have
\[
 f\in S_m^+
 \quad\Longleftrightarrow\quad q_{m,p}(f)=1
 \quad\Longleftrightarrow\quad \|J_m f\|_{C(Y_m)}=1.
\]
Since $J_m(K_m)=\mathcal A_m^+$, the map $J_m$ carries $S_m^+$ onto the positive supremum-norm unit sphere $\Sigma_m^+$.
Distances on this sphere are measured by $d_m$, so $J_m$ is an isometry from $(S_m^+,\rho_{m,p})$ onto $(\Sigma_m^+,d_m)$.

For $y\in Y_m$, let $\delta_y$ denote the point-evaluation functional defined by $\delta_y(F)=F(y)$. 
We define the \textit{peak fiber} and \textit{zero fiber} at \(y\), respectively, by
\[
 \mathcal P_y:=\{F\in\Sigma_m^+:F(y)=1\},
 \qquad
 \mathcal Z_y:=\{F\in\Sigma_m^+:F(y)=0\}.
\]

For $m=0$ the coordinate space is simply $\mathcal A_0=C^\infty([0,1],\R)$ with $\Sigma_0^+=S_0^+$.
In each case $d_m(F,G)$ is the maximum between $\|F-G\|_{C(Y_m)}$ and $H_p(J_m^{-1}(F-G))$.
\begin{example}\label{ex:topology}
The topology $\mathcal{T}_{\rp}$ on $\Sp$ is strictly finer than the topology $\mathcal{T}_\infty$ on $\Sp$.
In particular, let $\omega_k=2\pi k$, and consider
\[
 g_k(x)=1-\frac{1}{2\sqrt{\omega_k}}(1+\sin(\omega_kx)).
\]
Then $g_k\in\Sp$ and $\supnorm{\one-g_k}=1/\sqrt{\omega_k}\to0$, whereas $q_{1,p}(\one-g_k)\ge\supnorm{g_k'}=\sqrt{\frac{\pi k}{2}}\to\infty$.
Also 
\[
\frac12\phi(q_{1,p}(\one-g_k))=\frac12\frac{q_{1,p}(\one - g_k)}{1+q_{1,p}(\one-g_k)}\longrightarrow\frac12. 
\]
Since this is one of the terms defining $\Hp$, Lemma~\ref{lem:tail} gives
\[
 \frac12\phi(q_{1,p}(\one-g_k))
 \le \Hp(\one-g_k)<\frac12.
\]
By the squeeze theorem, $\Hp(\one-g_k)\to1/2$.
Since $q_0(\one-g_k)=\supnorm{\one-g_k}=1/\sqrt{\omega_k}\to0$, we conclude that $\rp(\one,g_k)\to1/2$.
\end{example}

The reconstruction of peak fibers in Lemma~\ref{lem:finite-midpoint-fibers} follows the method of Leung, Ng, and Wong~\cite{LeungNgWong}.
We verify that the method applies to $\Sigma_m^+$ equipped with $d_m$.

First, the bound $H_p<1/2$ gives
\[
 d_m(F,G)<1
 \quad\Longleftrightarrow\quad
 \|F-G\|_{C(Y_m)}<1
 \qquad(F,G\in\Sigma_m^+).
\]
Thus the distance-one relation is the same in both metrics.
This criterion also appears in~\cite[Lemma~2.1]{EzumiLinMiura}.
Second, the separating functions are chosen smooth on the interval component, with independently prescribed values at the isolated points.
This keeps all test functions in $\Sigma_m^+$.
Third, we pass from strictly positive members of each peak fiber to arbitrary members by using $(F+\varepsilon)/(1+\varepsilon)$.
The inverse coordinate images of these approximants converge in the smooth topology as $\varepsilon\downarrow0$, so the approximants converge to $F$ in $d_m$.
Since $\Phi$ preserves $d_m$ and $\|U-V\|_{C(Y_m)}\le d_m(U,V)$, their images converge uniformly to $\Phi F$.
This convergence preserves the prescribed peak value.

Once the peak and zero fibers have been identified, Lemma~\ref{lem:finite-midpoint-evaluations} recovers point evaluations from $d_m$-distances to these fibers.
The resulting composition representation then allows us to determine which coordinate maps preserve the ambient metric.

\begin{lemma}\label{lem:finite-midpoint-fibers}
If $\Phi:\Sigma_m^+\to\Sigma_m^+$ is a surjective $d_m$-isometry, there is a unique bijection $\tau:Y_m\to Y_m$ such that
\[
 \Phi(\mathcal P_y)=\mathcal P_{\tau(y)},
 \qquad
 \Phi(\mathcal Z_y)=\mathcal Z_{\tau(y)}
 \quad(y\in Y_m).
\]
\end{lemma}
\begin{proof}
We follow the reconstruction argument of
\cite[Lemma~14 and Theorem~15]{LeungNgWong}.
For $F\in\Sigma_m^+$, write
\[
 M(F)=\{y\in Y_m:F(y)=1\},\qquad
 Z(F)=\{y\in Y_m:F(y)=0\}.
\]
Both sets are compact, and $M(F)$ is nonempty.
Since $H_p<1/2$,
for every $F,G\in\Sigma_m^+$ we have
\begin{equation}\label{eq:jet-distance-one}
\begin{aligned}
 d_m(F,G)<1
 \quad\Longleftrightarrow\quad \|F-G\|_{C(Y_m)}<1
 \quad\Longleftrightarrow\quad
 M(F)\cap Z(G)=Z(F)\cap M(G)=\varnothing.
\end{aligned}
\end{equation}
The second equivalence follows from compactness of $Y_m$, since $F$ and $G$ take values in $[0,1]$.

We shall use the following separating functions on $Y_m$.
Given a nonempty compact set $K\subseteq Y_m$ and a point $t\notin K$, we may choose $V\in\Sigma_m^+$ such that $V=1$ on $K$ and $V(t)=0$.
In fact, we may use a smooth cutoff on the interval component and prescribe the values at the isolated points independently.
For each $y\in Y_m$, we also choose a strictly positive $U_y\in\Sigma_m^+$ whose unique peak is $y$.
If $y$ is isolated, take $U_y(y)=1$ and $U_y=1/2$ elsewhere.
If $y\in[0,1]$, put
\[
 U_y(t)=1-\tfrac12(t-y)^2\qquad(t\in[0,1])
\]
and assign the value $1/2$ at each isolated point.
Thus $M(U_y)=\{y\}$ and $Z(U_y)=\varnothing$ in either case.

For every $\mathcal E\subseteq\Sigma_m^+$, put
\[
 \mathcal E^\#
 :=\{H\in\Sigma_m^+:d_m(H,G)<1
                   \text{ for every }G\in\mathcal E\},
 \qquad
 \mathcal E^{\#\#}=(\mathcal E^\#)^\#.
\]
Since $\Phi$ is a surjective isometry,
\begin{equation}\label{eq:jet-sharp-preservation}
 \Phi(\mathcal E^\#)=(\Phi(\mathcal E))^\#.
\end{equation}
Let $\mathcal Q=\{F\in\Sigma_m^+:Z(F)=\varnothing\}$.
We first show that $F\in\mathcal Q$ if and only if
\begin{equation}\label{eq:jet-strict-positive}
 \{F\}^\#\cap\mathcal E^\#\ne\varnothing
 \quad\text{for every nonempty }\mathcal E\subseteq\Sigma_m^+
 \text{ with }\mathcal E^\#\ne\varnothing.
\end{equation}

Suppose that $F\in\mathcal Q$ and choose $H\in\mathcal E^\#$.
Let $\mathbf 1_{Y_m}$ denote the constant function with value one on $Y_m$, and put $H_0=(\mathbf 1_{Y_m}+H)/2$.
Since $M(H_0)=M(H)$ and $Z(H_0)=\varnothing$, criterion~\eqref{eq:jet-distance-one} gives $H_0\in\mathcal E^\#$.
Also, $F$ and $H_0$ are strictly positive, so $d_m(F,H_0)<1$.
This proves~\eqref{eq:jet-strict-positive}.
Conversely, suppose that $Z(F)\ne\varnothing$.
For each $t\in Y_m\setminus Z(F)$, choose $V_t\in\Sigma_m^+$ such that $V_t=1$ on $Z(F)$ and $V_t(t)=0$.
Put
\[
 \mathcal E=\{V_t:t\in Y_m\setminus Z(F)\}.
\]
Since $F$ attains the value one, $\mathcal E$ is nonempty.
The function $H_0=\mathbf 1_{Y_m}-F/2$ belongs to $\mathcal E^\#$, because $M(H_0)=Z(F)$ and $Z(H_0)=\varnothing$.
On the other hand, if $H\in\mathcal E^\#$, then $H(t)<1$ for every $t\notin Z(F)$, since $V_t(t)=0$.
Thus $M(H)\subseteq Z(F)$.
As $M(H)$ is nonempty, criterion~\eqref{eq:jet-distance-one} gives $d_m(F,H)=1$.
Hence $\{F\}^\#\cap\mathcal E^\#=\varnothing$.
This proves the claimed equivalence.
By~\eqref{eq:jet-sharp-preservation}, it follows that $\Phi(\mathcal Q)=\mathcal Q$.

Next we show that, for every $G\in\mathcal Q$,
\begin{equation}\label{eq:jet-double-sharp}
 \{G\}^{\#\#}
 =\{F\in\mathcal Q:M(F)\subseteq M(G)\}.
\end{equation}
Let $F\in\mathcal Q$ satisfy $M(F)\subseteq M(G)$.
Every $H\in\{G\}^\#$ is strictly positive on $M(G)$, and hence on $M(F)$.
Since $Z(F)=\varnothing$, it follows from~\eqref{eq:jet-distance-one} that $d_m(F,H)<1$.
Thus $F\in\{G\}^{\#\#}$.

Conversely, let $F\in\{G\}^{\#\#}$. If $t\in Z(F)$, the strictly positive function $U_t$ belongs to $\{G\}^\#$, whereas $U_t(t)=1$ and $F(t)=0$ give $d_m(F,U_t)=1$.
This contradiction shows that $F\in\mathcal Q$. If $t\in M(F)\setminus M(G)$, choose $V\in\Sigma_m^+$ equal to one on $M(G)$ and zero at $t$.
Then $V\in\{G\}^\#$, since $G$ is strictly positive, but $d_m(F,V)=1$.
This again contradicts $F\in\{G\}^{\#\#}$.
Hence $M(F)\subseteq M(G)$, proving~\eqref{eq:jet-double-sharp}.

For $F,G\in\mathcal Q$, equation~\eqref{eq:jet-double-sharp} gives
\[
 M(F)\subseteq M(G)
 \quad\Longleftrightarrow\quad
 \{F\}^{\#\#}\subseteq\{G\}^{\#\#}.
\]
Indeed, the forward implication follows by inclusion of peak sets.
For the reverse implication, use $F\in\{F\}^{\#\#}$.
Applying~\eqref{eq:jet-sharp-preservation} twice and using $\Phi(\mathcal Q)=\mathcal Q$, we obtain
\begin{equation}\label{eq:jet-peak-order}
 M(F)\subseteq M(G)
 \quad\Longleftrightarrow\quad
 M(\Phi F)\subseteq M(\Phi G)
 \qquad(F,G\in\mathcal Q).
\end{equation}

For $y\in Y_m$, put $\mathcal Q_y=\mathcal Q\cap\mathcal P_y$.
Fix $y$ and let $F_1,\ldots,F_n\in\mathcal Q_y$.
Their average $H=n^{-1}\sum_{j=1}^n F_j$ belongs to $\mathcal Q_y$ and satisfies
\[
 M(H)=\bigcap_{j=1}^n M(F_j).
\]
By~\eqref{eq:jet-peak-order},
\[
 \varnothing\ne M(\Phi H)
 \subseteq\bigcap_{j=1}^n M(\Phi F_j).
\]
Thus the compact sets $M(\Phi F)$, $F\in\mathcal Q_y$, have the finite intersection property.
Compactness of $Y_m$ gives a point $z$ in their common intersection.
Since $\Phi$ preserves $\mathcal Q$, we have $\Phi(\mathcal Q_y)\subseteq\mathcal Q_z$.
Applying the same argument to $\Phi^{-1}$ at $z$, we find $w\in Y_m$ such that
\[
 \mathcal Q_y
 \subseteq\Phi^{-1}(\mathcal Q_z)
 \subseteq\mathcal Q_w.
\]
Now $U_y\in\mathcal Q_y$ and $M(U_y)=\{y\}$, so $w=y$.
Both inclusions are therefore equalities, and $\Phi(\mathcal Q_y)=\mathcal Q_z$.

Distinct points have distinct fibers $\mathcal Q_y$, as witnessed by the functions $U_y$.
Hence $z$ is uniquely determined by $y$.
Writing $z=\tau(y)$, we obtain an injective map $\tau:Y_m\to Y_m$.
The same construction for $\Phi^{-1}$ shows that every $\mathcal Q_z$ is the image of some $\mathcal Q_y$.
Thus $\tau$ is bijective and
\[
 \Phi(\mathcal Q_y)=\mathcal Q_{\tau(y)}
 \qquad(y\in Y_m).
\]

We extend this identity to all peak fibers.
Let $F\in\mathcal P_y$ and, for $\varepsilon>0$, put
\[
 F_\varepsilon
 =\frac{F+\varepsilon\mathbf 1_{Y_m}}{1+\varepsilon}.
\]
Then $F_\varepsilon\in\mathcal Q_y$. Moreover,
\[
 J_m^{-1}(F_\varepsilon-F)
 =\frac{\varepsilon}{1+\varepsilon}
   J_m^{-1}(\mathbf 1_{Y_m}-F)
 \longrightarrow0
\]
in the smooth topology as $\varepsilon\downarrow0$.
By Theorem~\ref{thm:metric}, $d_m(F_\varepsilon,F)\to0$.
Hence
\[
 \|\Phi F_\varepsilon-\Phi F\|_{C(Y_m)}
 \le d_m(\Phi F_\varepsilon,\Phi F)
 =d_m(F_\varepsilon,F)\longrightarrow0.
\]
Since $(\Phi F_\varepsilon)(\tau(y))=1$ for every $\varepsilon>0$, passing to the limit gives $(\Phi F)(\tau(y))=1$.
Thus $\Phi(\mathcal P_y)\subseteq\mathcal P_{\tau(y)}$.
Applying the same argument to $\Phi^{-1}$ gives the reverse inclusion.

Finally, criterion~\eqref{eq:jet-distance-one} gives
\[
 F(y)=0
 \quad\Longleftrightarrow\quad
 \{F\}^\#\cap\mathcal P_y=\varnothing.
\]
To see this, note that if $F(y)=0$, every member of $\mathcal P_y$ has distance one from $F$.
If $F(y)>0$, then $U_y\in\{F\}^\#\cap\mathcal P_y$.
By~\eqref{eq:jet-sharp-preservation} and the peak-fiber identity,
\[
 F(y)=0\quad\Longleftrightarrow\quad(\Phi F)(\tau(y))=0.
\]
Surjectivity of $\Phi$ now gives $\Phi(\mathcal Z_y)=\mathcal Z_{\tau(y)}$ for every $y\in Y_m$.
Since the peak fibers at distinct points are distinct, their images determine $\tau$ uniquely.
\end{proof}

At order zero, the following observation explains the use of smooth approximants in the evaluation formula.
The final step in the proof of~\cite[Theorem~2.1]{EzumiLinMiura} cannot be transferred directly to the present setting.
In particular, the smooth analogue of~\cite[Lemma~2.8]{EzumiLinMiura} fails.
For example, take $v(t)=t$ and $y_0=1/2$.
If a smooth positive unit vector $u$ satisfies $u(y_0)=1$, then $u'(y_0)=0$.
Consequently,
\[
\left(v+\frac12u\right)'(y_0)=1,
\]
so $v+\frac12u$ cannot belong to $\Sp$, since it takes the value one at the interior point $y_0$.

\begin{lemma}\label{lem:finite-midpoint-evaluations}
For $F\in\Sigma_m^+$ and $y\in Y_m$, set
\[
 \mathfrak p_y(F)=\inf_{G\in\mathcal P_y}d_m(F,G),
 \qquad
 \mathfrak z_y(F)=\inf_{G\in\mathcal Z_y}d_m(F,G).
\]
Then
\begin{equation}\label{eq:finite-midpoint-evaluation}
 F(y)=
 \begin{cases}
 1-\mathfrak p_y(F),&\mathfrak p_y(F)\ge1/2,\\
 \mathfrak z_y(F),&\mathfrak p_y(F)<1/2.
 \end{cases}
\end{equation}
Consequently, if $\Phi$ and $\tau$ are as in Lemma~\ref{lem:finite-midpoint-fibers}, then
\begin{equation}\label{eq:finite-midpoint-point-action}
 (\Phi F)(\tau(y))=F(y)
 \qquad(F\in\Sigma_m^+,\ y\in Y_m).
\end{equation}
\end{lemma}
\begin{proof}
First, for every $y\in Y_m$ and every open neighborhood $U$ of $y$, there exists $\eta\in\mathcal A_m$ such that
\[
 0\le\eta\le1,\qquad \eta(y)=1,\qquad
 \supp\eta\subset U.
\]
If $y$ is isolated, take $\eta$ to be the indicator of $\{y\}$.
If $y$ belongs to the interval component, choose a smooth cutoff there, supported in $U\cap[0,1]$, and assign zero at every isolated point.
Consider the two smooth test functions
\[
 G=F+(1-F)\eta
 \qquad\text{and}\qquad
 G=F(1-\eta)
\]
The first construction belongs to $\mathcal P_y$ when $\eta$ is supported where $F>F(y)-\varepsilon$.
For the second construction, when $0<F(y)<1$, choose $\eta$ supported where $F<F(y)+\varepsilon$ and away from $M(F)$.
Then $G=F(1-\eta)$ belongs to $\mathcal Z_y$ and remains one on $M(F)$.
If $F(y)=0$, take $G=F$; if $F(y)=1$, use any member of $\mathcal Z_y$.
Both fibers are nonempty: $\one_{Y_m}\in\mathcal P_y$, and a cutoff supported near a point different from $y$ gives a member of $\mathcal Z_y$.
Evaluation at $y$ gives the lower bounds below.
The two cutoffs give the corresponding upper bounds within $\varepsilon$, since the changes are $(1-F)\eta$ and $F\eta$ on their respective supports.
Letting $\varepsilon\downarrow0$ gives
\[
 \inf_{G\in\mathcal P_y}\|F-G\|_{C(Y_m)}
 =1-F(y),\qquad
 \inf_{G\in\mathcal Z_y}\|F-G\|_{C(Y_m)}
 =F(y).
\]

Second, the definition of $d_m$ gives
\[
 d_m(F,G)
 =\max\left\{
   \|F-G\|_{C(Y_m)},
   H_p\bigl(J_m^{-1}(F-G)\bigr)
 \right\}.
\]
Since $H_p(h)<1/2$ for every smooth $h$, we have
\[
 \begin{aligned}
 \|F-G\|_{C(Y_m)}\ge1/2
 &\ \Longrightarrow\
 d_m(F,G)=\|F-G\|_{C(Y_m)},\\
 \|F-G\|_{C(Y_m)}<1/2
 &\ \Longrightarrow\
 d_m(F,G)<1/2.
 \end{aligned}
\]

If $1-F(y)\ge1/2$, every member of $\mathcal P_y$ has uniform distance at least $1/2$ from $F$, hence $\mathfrak p_y(F)=1-F(y)$. Otherwise the uniform infimum gives a member at distance strictly below $1/2$ in both metrics.
The same argument for $\mathcal Z_y$ gives $\mathfrak z_y(F)=F(y)$ whenever $F(y)\ge1/2$.
These two cases prove the evaluation formula.
Finally, the isometry property of $\Phi$ and the fiber identities in Lemma~\ref{lem:finite-midpoint-fibers} imply
\[
 \mathfrak p_{\tau(y)}(\Phi F)=\mathfrak p_y(F),
 \qquad
 \mathfrak z_{\tau(y)}(\Phi F)=\mathfrak z_y(F).
\]
Applying \eqref{eq:finite-midpoint-evaluation} to $F$ at $y$ and to $\Phi F$ at $\tau(y)$ gives \eqref{eq:finite-midpoint-point-action}.
\end{proof}

\begin{proposition}\label{prop:finite-midpoint-composition}
Suppose $m\ge 1$ and $p=\infty$.
Every surjective $d_m$-isometry $\Phi:\Sigma_m^+\to\Sigma_m^+$ has the form
\begin{equation}\label{eq:finite-midpoint-coordinate-form}
 \Phi(b_0,\ldots,b_{m-1},g)
 =\bigl(b_{\pi(0)},\ldots,b_{\pi(m-1)},
        g\circ\sigma\bigr),
\end{equation}
where $\pi$ is a permutation of $\{0,\ldots,m-1\}$ and $\sigma:[0,1]\to[0,1]$ is a smooth diffeomorphism.
\end{proposition}
\begin{proof}
By Lemma~\ref{lem:finite-midpoint-evaluations}, there is a bijection $\tau:Y_m\to Y_m$ such that
\[
 (\Phi F)(\tau(y))=F(y)
 \qquad(F\in\Sigma_m^+,\ y\in Y_m).
\]
Put $\theta=\tau^{-1}$.
Then $\Phi F=F\circ\theta$ for each $F\in\Sigma_m^+$.

We first show that composition with $\theta$ and $\theta^{-1}$ preserves $\mathcal A_m$.
By Proposition~\ref{prop:finite-jet-span} and the identification $J_m$, the set $\Sigma_m^+$ spans $\mathcal A_m$ over $\R$.
Thus, for every $F\in\mathcal A_m$, there are $c_1,\ldots,c_N\in\R$ and $F_1,\ldots,F_N\in\Sigma_m^+$ such that $F=\sum_{j=1}^N c_jF_j$. 
Consequently,
\[
 F\circ\theta
 =\sum_{j=1}^N c_j(F_j\circ\theta)
 =\sum_{j=1}^N c_j\Phi F_j
 \in\mathcal A_m.
\]
Moreover, surjectivity of $\Phi$ gives
\[
 \Phi^{-1}G=G\circ\theta^{-1}
 \qquad(G\in\Sigma_m^+).
\]
The same spanning argument therefore shows that $F\circ\theta^{-1}\in\mathcal A_m$ for every $F\in\mathcal A_m$.

We next prove that $\theta$ is a homeomorphism.
The space $Y_m$ is compact and metrizable, and $\mathcal A_m$ separates its points.
In particular, the characteristic functions of the isolated singletons distinguish those points from all other points, while the function that vanishes at the isolated points and equals $t\mapsto t$ on the interval distinguishes points of the interval.

Let $y_k\to y$ in $Y_m$.
If $\theta(y_k)$ did not converge to $\theta(y)$, compactness and metrizability would yield a subsequence $\theta(y_{k_\ell})$ converging to some $z\ne\theta(y)$.
For every $F\in\mathcal A_m$, both $F$ and $F\circ\theta$ are continuous, so
\[
 F(z)
 =\lim_{\ell\to\infty}F(\theta(y_{k_\ell}))
 =F(\theta(y)).
\]
This contradicts point separation by $\mathcal A_m$.
Hence $\theta$ is continuous.
Applying the same argument to $\theta^{-1}$ proves that $\theta$ is a homeomorphism.

A homeomorphism preserves isolated points.
Thus $\theta$ permutes $\{a_0,\ldots,a_{m-1}\}$ and maps the interval component onto itself.
There are therefore a permutation $\pi$ of $\{0,\ldots,m-1\}$ and a homeomorphism $\sigma:[0,1]\to[0,1]$ such that
\[
 \theta(a_j)=a_{\pi(j)}
 \quad(0\le j<m),
 \qquad
 \theta|_{[0,1]}=\sigma.
\]

It remains to verify smoothness of $\sigma$ and its inverse.
Let $I\in\Sigma_m^+$ be the function that vanishes at all isolated points and satisfies $I(t)=t$ on the interval.
The interval component of $\Phi I=I\circ\theta$ is $\sigma$.
Since $\Phi I\in\mathcal A_m$, the function $\sigma$ is smooth.
Similarly, the interval component of $\Phi^{-1}I=I\circ\theta^{-1}$ is $\sigma^{-1}$, so $\sigma^{-1}$ is smooth.
Hence $\sigma$ is a smooth diffeomorphism of $[0,1]$.

Finally, if $F=(b_0,\ldots,b_{m-1},g)\in\Sigma_m^+$, then
\[
 (\Phi F)(a_j)=F(\theta(a_j))=b_{\pi(j)}
 \qquad(0\le j<m),
\]
and the interval component of $\Phi F$ is $g\circ\sigma$.
Therefore
\[
 \Phi(b_0,\ldots,b_{m-1},g)
 =\bigl(b_{\pi(0)},\ldots,b_{\pi(m-1)},g\circ\sigma\bigr),
\]
as asserted.
\end{proof}

We now identify the admissible interval maps at order zero.
Let $\Delta^+:S_0^+\to S_0^+$ be a surjective $\rho_{0,p}$-isometry.
Since $J_0$ is the identity and $Y_0=[0,1]$, Lemma~\ref{lem:finite-midpoint-evaluations} gives a bijection $\tau:[0,1]\to[0,1]$ such that
\[
(\Delta^+ f)(\tau(t))=f(t)
\qquad(f\in S_0^+,\ t\in[0,1]).
\]
Putting $\sigma=\tau^{-1}$, we obtain
\[
\Delta^+ f=f\circ\sigma
\qquad(f\in S_0^+).
\]
Moreover, the function $u(t)=t$ belongs to $S_0^+$, and hence
\[
\sigma=\Delta^+u\in C^\infty([0,1],\R),
\qquad
\sigma^{-1}=(\Delta^+)^{-1}u
\in C^\infty([0,1],\R).
\]
Thus $\sigma$ is a smooth diffeomorphism of $[0,1]$.

Recall that
\[
 q_{1,p}(h)=\begin{cases}
 (|h(x_0)|^p+\supnorm{h'}^p)^{1/p},\quad& p<\infty,\\
 \max\{|h(x_0)|,\supnorm{h'}\},\quad&p=\infty.
 \end{cases}
\]
We introduce the class
\[
 \mathcal H=\{h\in  C^\infty([0,1], \R):h\ge0,\ \min h=0,\ q_{1,p}(h)>6\supnorm h\}.
\]

\begin{lemma}\label{lem:isolation}
Suppose that \(\Delta^+:\Sp\to\Sp\) is a surjective \(\rp\)-isometry of the form \(\Delta^+(f)=f\circ\sigma\).
Then
\begin{equation}\label{eq:firstnorm}
 q_{1,p}(h\circ\sigma)=q_{1,p}(h)\qquad(h\in\mathcal H).
\end{equation}
\end{lemma}
\begin{proof}
Fix $h\in\mathcal H$ and put $s=2/q_{1,p}(h)$.
Then $s\supnorm h<1/3$, and $g=\one-sh\in\Sp$.
The first derivative contribution to $\rp(sh,0)$ equals
\[
 \frac12\frac{sq_{1,p}(h)}{1+sq_{1,p}(h)}=\frac13.
\]
All terms with $n\ge2$ are at most $1/4$.
The constant $6$ in the definition of $\mathcal H$ ensures that the supremum-norm term is smaller than $1/3$, so the first contribution is isolated and $\rp(\one,g)=1/3$.

By the composition formula, $\Delta^+(\one)=\one$ and $\Delta^+ (g)=\one-s(h\circ\sigma)$.
Hence $\rp(s(h\circ\sigma),0)=1/3$.
Composition by a homeomorphism preserves the supremum norm, and all terms with $n\ge2$ are still at most $1/4$.
Therefore
\[
 \frac12\frac{sq_{1,p}(h\circ\sigma)}{1+sq_{1,p}(h\circ\sigma)}=\frac13,
\]
which gives $sq_{1,p}(h\circ\sigma)=2$ and proves~\eqref{eq:firstnorm}.
\end{proof}

\begin{lemma}\label{lem:oscillation}
Let $\sigma:[0,1]\to[0,1]$ be a $C^1$ diffeomorphism and set
\[
 M_k=\sup_{x\in[0,1]}\abs{\sin(2\pi k\sigma(x))\sigma'(x)}.
\]
Then $M_k\to\supnorm{\sigma'}$ as $k\to\infty$.
\end{lemma}
\begin{proof}
The inequality $M_k\le\supnorm{\sigma'}$ is immediate.
Fix $\varepsilon>0$.
By continuity of $\sigma'$, there is a nondegenerate interval $U\subset[0,1]$ such that $|\sigma'(x)|>\supnorm{\sigma'}-\varepsilon$ for every $x\in U$.
Since $\sigma$ is a diffeomorphism, $\sigma(U)$ contains a nondegenerate interval $V$.

The points at which $|\sin(2\pi ky)|=1$ are spaced by $1/(2k)$.
Consequently, for every sufficiently large integer $k$, there exists $y_k\in V$ such that $|\sin(2\pi ky_k)|=1$.
Setting $x_k=\sigma^{-1}(y_k)\in U$, we obtain
\[
 M_k\ge|\sigma'(x_k)\sin(2\pi ky_k)|=|\sigma'(x_k)|>\supnorm{\sigma'}-\varepsilon.
\]
Since $\varepsilon>0$ is arbitrary, the conclusion follows.
\end{proof}

\begin{lemma}\label{lem:separatingH}
If $a,b\in[0,1]$ and $a\ne b$, there exists $h\in\mathcal H$ with $h(a)\ne h(b)$.
\end{lemma}
\begin{proof}
Choose $c\in(0,1)\setminus\{a,b,(a+b)/2\}$ and put $b_0(x)=(x-c)^2$.
Then $b_0(c)=0$ and $b_0(a)\ne b_0(b)$.
Choose a nonzero smooth function $\chi\ge0$ supported in an open subinterval of $(0,1)$ disjoint from $\{a,b,c\}$, and set
\[
 h_k(x)=b_0(x)+\chi(x)(1+\sin(2\pi kx)).
\]
These functions are nonnegative, vanish at $c$, and satisfy $h_k(a)=b_0(a)\ne b_0(b)=h_k(b)$.
Moreover, $\supnorm{h_k}\le\supnorm{b_0}+2\supnorm\chi$ uniformly in $k$.
On a fixed nondegenerate closed subinterval $J$, we have $\chi\ge\delta>0$.
For large $k$, choose $x_k\in J$ with $|\cos(2\pi kx_k)|=1$.
Differentiation gives
\[
 |h_k'(x_k)|\ge2\pi k\delta-\supnorm{b_0'}-2\supnorm{\chi'}.
\]
Thus $\supnorm{h_k'}\to\infty$, and $h_k\in\mathcal H$ for all sufficiently large $k$.
\end{proof}

Let $R$ denote the reflection operator on $C^\infty([0,1],\F)$, defined by $Rh:=h\circ r$.

\begin{proposition}\label{prop:finite}
Suppose $1\le p<\infty$.
The restriction $R|_{\Sp}$ is a surjective $\rp$-isometry if and only if $x_0=1/2$.
When $x_0=1/2$, the operator $R$ is a surjective $\F$-linear $\rp$-isometry on $ C^\infty([0,1], \F)$.
\end{proposition}
\begin{proof}
If $R|_{\Sp}$ is an isometry, Lemma~\ref{lem:isolation} and invariance of $\supnorm{h'}$ under reflection give
\[
 |h(1-x_0)|^p+\supnorm{h'}^p=|h(x_0)|^p+\supnorm{h'}^p
 \qquad(h\in\mathcal H).
\]
Since $h\ge0$, we have $h(1-x_0)=h(x_0)$ for all such $h$.
Lemma~\ref{lem:separatingH} forces $1-x_0=x_0$.
Conversely, at the centered anchor,
\[
 (h\circ r)^{(j)}(1/2)=(-1)^jh^{(j)}(1/2),\qquad
 \supnorm{(h\circ r)^{(n)}}=\supnorm{h^{(n)}}.
\]
Thus $R$ preserves every $\qn$ and $q_0$, and therefore $\rp$.
It is $\F$-linear and satisfies $R^2=\id$.
\end{proof}

\begin{proposition}\label{prop:anchorfree}
For $h\in C^\infty([0,1], \F)$, 
\begin{equation}\label{eq:anchorfree}
 \rho_\infty(h,0)=\max\left\{\supnorm{h},\,
                    \sup_{n\ge1}2^{-n}\frac{\supnorm{h^{(n)}}}{1+\supnorm{h^{(n)}}}\right\}.
\end{equation}
Hence $\rho_\infty$ is independent of the anchor, and reflection is a surjective $\F$-linear $\rho_\infty$-isometry for every anchor.
\end{proposition}
\begin{proof}
The inequality $q_{n,\infty}(h)\ge \supnorm{h^{(n)}}$ gives one direction.
For the other, each anchored coordinate satisfies $|h^{(j)}(x_0)|\le \supnorm{h^{(j)}}$.
If $j=0$, then $2^{-n}\phi(\supnorm{h})\le \supnorm{h}$.
If $1\le j<n$, then $2^{-n}\phi(\supnorm{h^{(j)}})\le2^{-j}\phi(\supnorm{h^{(j)}})$.
Since $\phi$ is increasing and preserves finite maxima, every contribution from $q_{n,\infty}$ is bounded by the right-hand side of~\eqref{eq:anchorfree}.
Reflection preserves every $\supnorm{h^{(n)}}$, proving the remaining assertions.
\end{proof}

We have defined $B_j(L_{x_0})$ and $L_{x_0}$ in \eqref{eq:L_and_B}.
For $m\ge1$ and either scalar field $\F\in\{\R,\C\}$, define
\[
 \mathcal D_m(g)
 :=\max\left\{\|g\|_\infty,
   \sup_{k\ge1}2^{-(m+k)}
      \phi(\|g^{(k)}\|_\infty)\right\}
 \qquad (g\in C^\infty([0,1],\F)).
\]
We establish the following lemma for both scalar fields.

\begin{lemma}\label{lem:finite-midpoint-splitting}
Let $m\ge1$, $x_0\in[0,1]$, and $\F\in\{\R,\C\}$.
Suppose that
\begin{equation}\label{eq:finite-anchor-condition}
 B_{m-1}(L_{x_0})\le2.
\end{equation}
For $f\in C^\infty([0,1],\F)$, write $J_mf=(b_0,\ldots,b_{m-1},g)$, where $b_j=f^{(j)}(x_0)$ and $g=f^{(m)}$.
Then
\begin{equation}\label{eq:finite-midpoint-split}
 \rho_{m,\infty}(f,0) = \max\left\{\max_{0\le j<m}|b_j|, \mathcal D_m(g)\right\}.
\end{equation}
\end{lemma}
\begin{proof}
Recall $q_{m,\infty}(f)=\max\left\{\max_{0\le j<m}|b_j|,\|g\|_\infty\right\}$.
For $1\le n<m$, Taylor's formula with integral remainder gives, for every $t\in[0,1]$,
\[
 f^{(n)}(t)
 =\sum_{j=n}^{m-1}
   \frac{f^{(j)}(x_0)}{(j-n)!}
   (t-x_0)^{j-n}
  +\int_{x_0}^{t}
   \frac{(t-s)^{m-n-1}}{(m-n-1)!}
   g(s)\,ds.
\]
Taking absolute values and using $|t-x_0|\le L_{x_0}$ gives, for either scalar field,
\[
 \|f^{(n)}\|_\infty
 \le q_{m,\infty}(f)
    \sum_{k=0}^{m-n}
    \frac{L_{x_0}^k}{k!}
 =B_{m-n}(L_{x_0})q_{m,\infty}(f).
\]
Since the anchored derivatives of orders below $n$ have absolute values at most $q_{m,\infty}(f)$, it follows that
\begin{equation}\label{eq:finite-anchor-lower-order}
 q_{n,\infty}(f)\le B_{m-n}(L_{x_0})q_{m,\infty}(f)
 \qquad(1\le n<m).
\end{equation}
For $n=1<m$, condition~\eqref{eq:finite-anchor-condition} and $\phi(u)\le u$ yield
\[
 2^{-1}\phi(q_{1,\infty}(f))
 \le\frac12B_{m-1}(L_{x_0})q_{m,\infty}(f)\le q_{m,\infty}(f).
\]
For $2\le n<m$, we use $L_{x_0}\le1$ and $B_{m-n}(L_{x_0})\le e^{L_{x_0}}\le e<4$ to obtain
\[
 2^{-n}\phi(q_{n,\infty}(f))
 \le2^{-n}B_{m-n}(L_{x_0})q_{m,\infty}(f)
 \le\frac e4 q_{m,\infty}(f)<q_{m,\infty}(f)\qquad(q_{m,\infty}(f)>0).
\]
The term with $n=m$ is at most $q_{m,\infty}(f)$.
Thus every term of $H_\infty(f)$ of order at most $m$ is bounded by the base norm $q_{m,\infty}(f)$.

For $k\ge1$, the anchored-coordinate definition gives
\[
 q_{m+k,\infty}(f)
 =\max\left\{
   \max_{0\le j<m}|b_j|,
   \max_{0\le j<k}|g^{(j)}(x_0)|,
   \|g^{(k)}\|_\infty
 \right\}.
\]
The terms $|b_j|$ and $|g(x_0)|$ are bounded by $q_{m,\infty}(f)$.
For $1\le j<k$, monotonicity of $\phi$ gives
\[
 2^{-(m+k)}\phi(|g^{(j)}(x_0)|)
 \le2^{-(m+j)}\phi(\|g^{(j)}\|_\infty)
 \le\mathcal D_m(g).
\]
The terminal derivative term is bounded by $\mathcal D_m(g)$ as well.
Hence 
\[
\rho_{m,\infty}(f,0)\le\max\{q_{m,\infty}(f),\mathcal D_m(g)\}.
\]
Conversely, $\rho_{m,\infty}(f,0)\ge q_{m,\infty}(f)$, and for each $k\ge1$ it is at least $2^{-(m+k)}\phi(\|g^{(k)}\|_\infty)$.
Since $q_{m,\infty}(f)\ge\|g\|_\infty$, the reverse inequality follows.
This proves~\eqref{eq:finite-midpoint-split}.
\end{proof}

\begin{corollary}\label{cor:symbol}
Under the hypotheses of Theorem~\ref{thm:main}, let $\sigma$ be the interval map induced by a surjective $\rho_{m,p}$-isometry $\Delta^+$ of $S_m^+$.
Then $\sigma=\id$ or $\sigma=r$, where $r(t)=1-t$.
\end{corollary}

\begin{proof}
For $m\ge1$, $\sigma$ is a smooth diffeomorphism by Proposition~\ref{prop:finite-midpoint-composition}.
For $m=0$, this follows from the order-zero composition representation established immediately before Lemma~\ref{lem:isolation}.
Put $h_k(t)=(1-\cos(2\pi kt))/(\pi k)$. Then $h_k(0)=0$, $0\le h_k\le 2/(\pi k)$, and $\supnorm{h_k'}=2$.

Suppose $p=\infty$.
Compare the sphere elements with zero finite coordinates and terminal components $\one$ and $\one-h_k$.
Lemma~\ref{lem:finite-midpoint-splitting} for $m\ge1$, and Proposition~\ref{prop:anchorfree} for $m=0$, give
\[
 \mathcal D_m(h_k)=\mathcal D_m(h_k\circ\sigma),
\]
where $\mathcal D_0$ denotes the expression in~\eqref{eq:anchorfree}.
Put $\alpha=2^{-(m+1)}$.
The first-derivative term on the left is $\alpha\phi(2)=2\alpha/3$.
On both sides the uniform term tends to zero and all higher-derivative terms are at most $\alpha/2$.
Hence $\supnorm{(h_k\circ\sigma)'}=2$ for all sufficiently large $k$.

If $p<\infty$, then $m=0$.
Since $q_{1,p}(h_k)\ge2>6\supnorm{h_k}$ for large $k$, Lemma~\ref{lem:isolation} gives
\[
 \supnorm{(h_k\circ\sigma)'}^p
 =2^p+|h_k(x_0)|^p-|h_k(\sigma(x_0))|^p
 \longrightarrow 2^p.
\]
Thus in both cases, by the chain rule and Lemma~\ref{lem:oscillation},
\[
 \supnorm{\sigma'}
 =\lim_{k\to\infty}\sup_{t\in[0,1]}
   |\sin(2\pi k\sigma(t))\,\sigma'(t)|
 =\lim_{k\to\infty}\tfrac12\supnorm{(h_k\circ\sigma)'}=1.
\]
Since $\sigma$ is a monotone bijection, $\int_0^1|\sigma'(t)|\,dt=|\sigma(1)-\sigma(0)|=1$.
Continuity gives $|\sigma'|\equiv1$, so $\sigma=\id$ or $\sigma=r$.
\end{proof}

\begin{proof}[\textbf{Proof of Theorem~\ref{thm:main}}]
Let $\Delta^+:S_m^+\to S_m^+$ be a surjective $\rho_{m,p}$-isometry and put $\Phi=J_m\Delta^+J_m^{-1}$.
Then $\Phi$ is a surjective $d_m$-isometry of $\Sigma_m^+$.
For $m\ge1$, Proposition~\ref{prop:finite-midpoint-composition} gives a permutation $\pi$ of $\{0,\ldots,m-1\}$ and an interval map $\sigma$ such that, whenever $J_mf=(b_0,\ldots,b_{m-1},g)$ and $f\in S_m^+$,
\[
 J_m(\Delta^+ f)=\bigl((b_{\pi(j)})_{0\le j<m},g\circ\sigma\bigr).
\]
For $m=0$, the same formula, with the finite block omitted, follows from Lemma~\ref{lem:finite-midpoint-evaluations}.
By Corollary~\ref{cor:symbol}, $\sigma=\id$ or $\sigma=r$.
If $p<\infty$, then $m=0$, and Proposition~\ref{prop:finite} shows that $\sigma=r$ is possible only when $x_0=1/2$.
Thus $\sigma$ is admissible.

Fix $\F\in\{\R,\C\}$. For any permutation $\pi$ and admissible $\sigma$, define $\widetilde{\Delta^+}$ by~\eqref{eq:finite-midpoint-extension}.
Since $J_m$ is a linear bijection, $\widetilde{\Delta^+}$ is $\F$-linear and has inverse given by $(\pi^{-1},\sigma^{-1})$.
Both maps preserve $S_m^+$, since permuting the finite coordinates and composing the interval coordinate with $\sigma$ preserve positivity and the maximum norm.

We verify that $\widetilde{\Delta^+}$ preserves the ambient metric.
For $m=0$, this follows from Propositions~\ref{prop:finite} and~\ref{prop:anchorfree}.
For $m\ge1$, we have $p=\infty$, and Lemma~\ref{lem:finite-midpoint-splitting} gives
\[
 \rho_{m,\infty}(f,0)
 =\max\left\{\max_{0\le j<m}|b_j|,\mathcal D_m(g)\right\},
 \qquad J_mf=(b_0,\ldots,b_{m-1},g).
\]
For $\sigma=\id$ or $r$, we have $\supnorm{(g\circ\sigma)^{(k)}}=\supnorm{g^{(k)}}$ for every $k\ge0$.
Hence $\mathcal D_m(g\circ\sigma)=\mathcal D_m(g)$, and the permutation preserves $\max_{0\le j<m}|b_j|$.
It follows that $\rho_{m,\infty}(\widetilde{\Delta^+}f,0)=\rho_{m,\infty}(f,0)$.
Linearity and translation invariance now give preservation of all distances.
This also proves the converse.

For the pair induced by $\Delta^+$, the coordinate formula gives $\widetilde{\Delta^+}|_{S_m^+}=\Delta^+$. Proposition~\ref{prop:finite-jet-span} shows that $S_m^+$ spans $C^\infty([0,1],\F)$ over $\F$, so the linear extension is unique.
Finally, its values on the finite coordinate vectors determine $\pi$, and its value on the vector with zero finite block and interval component $g(t)=t$ determines $\sigma$.
\end{proof}

\begin{corollary}[All anchors in the first two orders] \label{cor:finite-loworder-all-anchors}
Let $m\in\{1,2\}$, $p=\infty$, $x_0\in[0,1]$, and $\F\in\{\R,\C\}$.
Every surjective $\rho_{m,\infty}$-isometry $\Delta^+:S_m^+\to S_m^+$ admits a unique $\F$-linear surjective $\rho_{m,\infty}$-isometric extension to $C^\infty([0,1],\F)$, with the coordinate form \eqref{eq:finite-midpoint-extension}.
\end{corollary}
\begin{proof}
For $m=1$, $B_{m-1}(L_{x_0})=1$.
For $m=2$, $B_{m-1}(L_{x_0})=1+L_{x_0}\le2$.
Thus~\eqref{eq:finite-anchor-condition} holds for every anchor, and Theorem~\ref{thm:main} applies.
\end{proof}

More generally, the condition holds in every order if $1-\log2\le x_0\le\log2$, since $B_{m-1}(L_{x_0})\le e^{L_{x_0}}\le2$.
For $m\ge1$, the extension problem remains open here when $1\le p<\infty$, and when $p=\infty$ and $B_{m-1}(L_{x_0})>2$.

\begin{remark}
The scalar field matters for uniqueness.
Under the hypotheses of Theorem~\ref{thm:main}, let $\widetilde{\Delta^+}$ be the $\C$-linear extension and let $J$ denote pointwise conjugation.
Both $\widetilde{\Delta^+}$ and $J\circ\widetilde{\Delta^+}$ are $\R$-linear surjective metric isometries on the $\C$-valued space: conjugation preserves every defining norm and fixes $S_m^+$ pointwise.
The two operators are distinct by surjectivity of $\widetilde{\Delta^+}$.
\end{remark}

\subsection{A compatible convex-ball metric without the extension property}\label{subsec:fail_extension}
One may ask whether the positive-sphere extension property holds for every translation-invariant metric on $C^\infty([0,1],\R)$ that induces the usual Fr\'echet topology, is unbounded, has convex balls, and increases continuously and strictly along every nonzero positive ray.
The following construction provides a negative answer.
Its positive unit sphere is the same set as for $\rp$, and reflection is an isometry of that sphere, but it has no $\R$-linear isometric extension.
Throughout this subsection, functions are $\R$-valued and the anchor is $0$.

Write
\[
 \widehat q_n(h)=\max\left\{\max_{0\le j<n}|h^{(j)}(0)|,\,
                            \supnorm{h^{(n)}}\right\}
 \qquad(n\ge1),
\]
and define
\begin{equation}\label{eq:countermetric}
 d(f,g)=\max\left\{\supnorm{f-g},\,
       \sup_{n\ge1}2^{-n}\frac{3\widehat q_n(f-g)}{1+3\widehat q_n(f-g)}\right\}.
\end{equation}
Thus the only change from $\rho_\infty$ at this anchor is multiplication of each anchored norm by $3$ inside the bounded terms.

\begin{theorem}[Failure of positive-sphere extension]\label{thm:counterexample}
The metric $d$ in~\eqref{eq:countermetric} is translation invariant, unbounded, and compatible with the usual Fr\'echet topology $\mathcal{T}$ on $ C^\infty([0,1], \R)$.
Its balls are convex, and its positive-radius balls centered at zero are balanced and absorbing.
For every $h\ne0$, the function $t\mapsto d(th,0)$ is continuous, strictly increasing, and unbounded on $[0,\infty)$.
Its positive radius-one sphere is $\Sp$.
Nevertheless, the map
\[
 \Delta^+:\Sp\to\Sp,\qquad (\Delta^+ f)(x)=f(1-x),
\]
is a surjective $d$-isometry that has no $\R$-linear $d$-isometric extension to $ C^\infty([0,1], \R)$.
\end{theorem}
\begin{proof}
Since $\widehat q_n=q_{n,\infty}$ at the anchor $0$, the proof of Theorem~\ref{thm:metric}, applied to $3\widehat q_n$, shows that $d$ is a translation-invariant metric inducing $\mathcal T$.
It is unbounded because $d(h,0)\ge\supnorm h$.

The argument of Proposition~\ref{prop:large} shows that every $d$-ball is convex.
The balls centered at zero are balanced by homogeneity of $\widehat q_n$ and monotonicity of $\phi$, and absorbing by continuity of scalar multiplication.
Similarly, the proof of Proposition~\ref{prop:rays} gives the asserted radial properties.

Since
\[
\sup_{n\ge1}2^{-n}\phi(3\widehat q_n(h))
\le\max\left\{
\frac12\phi(3\widehat q_1(h)),\frac14
\right\}<\frac12,
\]
the positive unit sphere of $d$ is $\Sp$.
Now let $f,g\in\Sp$, put $h=f-g$, and write $a_j=\supnorm{h^{(j)}}$ for $j\ge0$.
Since both $f$ and $g$ are nonnegative and attain the value $1$, the function $h$ has a zero in $[0,1]$.
Hence $|h(0)|\le a_1$, and therefore
\[
a_n\le\widehat q_n(h)
\le\max_{1\le j\le n}a_j
\qquad(n\ge1).
\]
By monotonicity of $\phi$, this implies
\[
2^{-n}\phi(3a_n)
\le 2^{-n}\phi(3\widehat q_n(h))
\le\sup_{j\ge1}2^{-j}\phi(3a_j).
\]
Taking the supremum over $n$ yields
\begin{equation}\label{eq:countersphere}
d(f,g)=\max\left\{a_0,\sup_{n\ge1}2^{-n}\phi(3a_n)\right\}.
\end{equation}
Since reflection preserves every $a_j$, the map $\Delta^+$ is a surjective $d$-isometry of $\Sp$.

Finally, let $Rh=h\circ r$ and take
\[
h(x)=\frac{4+x}{40}.
\]
Then
\[
\supnorm h=\supnorm{Rh}=\frac18,
\qquad
\widehat q_n(h)=\frac1{10},
\qquad
\widehat q_n(Rh)=\frac18
\quad(n\ge1).
\]
Consequently,
\[
d(h,0)=\frac18,
\qquad
d(Rh,0)=\frac3{22}.
\]
Thus $R$ is not a global $d$-isometry.
Since $\Sp$ spans $C^\infty([0,1],\R)$ over $\R$ by Proposition~\ref{prop:finite-jet-span} with $m=0$, every $\R$-linear extension of $\Delta^+$ must equal $R$.
Hence no $\R$-linear $d$-isometric extension exists.
\end{proof}

\begin{remark}
The restriction~\eqref{eq:countersphere} is independent of anchored values because every difference of two positive sphere elements has a zero.
That fact need not hold for arbitrary functions in $ C^\infty([0,1], \R)$.
The extension fails because the ambient metric can detect anchored information that its restriction to $\Sp$ cannot detect.
For the original metric $\rho_\infty$, Proposition~\ref{prop:anchorfree} eliminates the anchored terms on the whole space, which is a stronger property.
The counterexample consequently separates uniqueness of a possible linear extension from existence of an isometric one.
It does not give a full-sphere isometry without an extension.
In fact, the two functions
\[
 f(x)=-\frac9{10}+\frac{19}{10}x,\qquad
 g(x)=-1+\frac{15}{8}x
\]
belong to the full real unit sphere and satisfy $f-g=(4+x)/40$.
The calculation in the proof gives $d(f,g)=1/8$ and $d(Rf,Rg)=3/22$, so reflection is not an isometry of that full sphere for the modified metric.
\end{remark}

\section{Full-sphere extensions and global isometries}\label{sec:full}
We equip the full uniform-norm unit sphere $S_{\F}$ with the order-zero metric $\rho_p=\rho_{0,p}$, where $\F\in\{\R,\C\}$.
To prove Theorem~\ref{thm:full}, we recover the pointwise action of each surjective isometry from distances on $S_{\F}$ and obtain its unique $\R$-linear isometric extension. We then use Corollary~\ref{cor:affine} to classify all surjective $\rho_p$-isometries of the ambient space.
Throughout the reconstruction argument, fix a surjective $\rp$-isometry $\Delta:S_{\F}\to S_{\F}$ and set
\[
 \mathbb T_{\F}=\{\lambda\in\F:|\lambda|=1\},\qquad
 F_{x,\lambda}:=\{f\in S_{\F}:f(x)=\lambda\},
\]
where we call $F_{x,\lambda}$ a \textit{phase fiber}.
Thus $\mathbb T_{\R}=\{-1,1\}$ and $\mathbb T_{\C}$ is the unit circle.
For $f\in S_{\F}$, define
\[
 \Pi(f):=\{(x,\lambda)\in[0,1]\times\mathbb T_{\F}:f(x)=\lambda\},
 \qquad D_2(f):=\{g\in S_{\F}:\rp(f,g)=2\}.
\]
The set $\Pi(f)$ is nonempty by compactness.
Here $\Pi(f)$ is called the \textit{boundary phase set}, and $D_2(f)$ is called the \textit{diametral set}.

\begin{lemma}\label{lem:boundaryfibers}
The minimal members of $\{D_2(f):f\in S_{\F}\}$, ordered by inclusion, are exactly $F_{x,\lambda}$.
Consequently, there are a bijection $\tau:[0,1]\to[0,1]$ and bijections $\alpha_x:\mathbb T_{\F}\to\mathbb T_{\F}$ such that
\begin{equation}\label{eq:boundarypermutation}
 \Delta(F_{x,\lambda})=F_{\tau(x),\alpha_x(\lambda)}.
\end{equation}
Each $\alpha_x$ preserves Euclidean distance.
For $\F=\R$ it has the form $\alpha_x(\lambda)=v_x\lambda$, where $v_x\in\{-1,1\}$.
For $\F=\C$ it has the form $v_x\lambda$ or $v_x\overline\lambda$, where $|v_x|=1$.
In either case, it extends to a $\R$-linear Euclidean isometry $\widetilde\alpha_x:\F\to\F$ fixing zero.
\end{lemma}
\begin{proof}
For $f,g\in S_{\F}$, we have $\supnorm f=\supnorm g=1$.
Since $\Hp(f-g)<1/2$,
\[
 \rp(f,g)=2
 \quad\Longleftrightarrow\quad
 \supnorm{f-g}=2.
\]
By compactness of $[0,1]$, the supremum norm is attained.
For $x\in[0,1]$, since $|f(x)|,|g(x)|\le1$, we have $|f(x)-g(x)|=2$ if and only if $f(x)=-g(x)\in\mathbb T_{\F}$.
Consequently,
\begin{equation}\label{eq:diametercriterion}
 \rp(f,g)=2
 \quad\Longleftrightarrow\quad
 f(x)=-g(x)\in\mathbb T_{\F}
 \text{ for some }x\in[0,1].
\end{equation}

By the definitions of $D_2(f)$, $\Pi(f)$, and $F_{x,\lambda}$, equation~\eqref{eq:diametercriterion} gives
\[
 D_2(f)
 =\bigcup_{(x,\lambda)\in\Pi(f)}
 F_{x,-\lambda}.
\]
We claim that
\begin{equation}\label{eq:diameterorder}
 D_2(f)\subseteq D_2(g)
 \quad\Longleftrightarrow\quad
 \Pi(f)\subseteq\Pi(g).
\end{equation}
The reverse implication follows from the union formula.
Conversely, suppose that $D_2(f)\subseteq D_2(g)$ and fix $(x,\lambda)\in\Pi(f)$.
Put
\[
 u_x(t)=1-\tfrac12(t-x)^2.
\]
Then $u_x$ has $x$ as its unique point of modulus one.
Since $f(x)=\lambda$, we have $-\lambda u_x\in D_2(f)$, and hence $-\lambda u_x\in D_2(g)$.
By \eqref{eq:diametercriterion}, there exists $y\in[0,1]$ such that
\[
 g(y)=\lambda u_x(y)
 \in\mathbb T_{\F}.
\]
Since $|u_x(y)|=1$ only when $y=x$, we obtain $g(x)=\lambda$.
Thus $(x,\lambda)\in\Pi(g)$, proving \eqref{eq:diameterorder}.

For every $x\in[0,1]$ and $\lambda\in\mathbb T_{\F}$, we have
\[
 \Pi(\lambda u_x)=\{(x,\lambda)\},
 \qquad
 D_2(\lambda u_x)=F_{x,-\lambda}.
\]
Hence every phase fiber occurs as a diametral set.

We next show that $D_2(f)$ is minimal under inclusion if and only if $\Pi(f)$ is a singleton.
Suppose that $\Pi(f)=\{(x,\lambda)\}$.
If $D_2(g)\subseteq D_2(f)$, then \eqref{eq:diameterorder} gives
\[
 \varnothing\ne\Pi(g)
 \subseteq\Pi(f)=\{(x,\lambda)\}.
\]
The nonemptiness follows from $\supnorm g=1$ and compactness of $[0,1]$.
Thus $\Pi(g)=\Pi(f)$, and \eqref{eq:diameterorder} yields $D_2(g)=D_2(f)$.
Therefore $D_2(f)$ is minimal.

Conversely, suppose that $\Pi(f)$ contains more than one element.
Choose $(x,\lambda)\in\Pi(f)$.
Then
\[
 \Pi(\lambda u_x)=\{(x,\lambda)\}
 \subsetneq\Pi(f).
\]
By both implications of \eqref{eq:diameterorder}, we obtain
\[
 D_2(\lambda u_x)\subsetneq D_2(f).
\]
Hence $D_2(f)$ is not minimal.
We conclude that the minimal diametral sets are precisely the phase fibers $F_{x,\lambda}$.

Since $\Delta$ is a surjective $\rp$-isometry,
\[
 \Delta(D_2(f))=D_2(\Delta f)
 \qquad(f\in S_{\F}).
\]
The same identity holds for $\Delta^{-1}$.
Therefore $\Delta$ permutes the minimal diametral sets, and hence the phase fibers.

We now recover their point labels.
We claim that two distinct phase fibers $F_{x,\lambda}$ and $F_{y,\mu}$ are disjoint if and only if $x=y$.
In particular, if $x=y$ and $\lambda\ne\mu$, no function can take both values at $x$.

Conversely, suppose that $x\ne y$.
Define
\[
 \eta(t):=
 \frac{(t-y)^2}
 {(t-x)^2+(t-y)^2},
 \qquad t\in[0,1].
\]
Then $\eta$ is smooth, $0\le\eta\le1$, and $\eta(x)=1$, $\eta(y)=0$.
Put $h=\lambda\eta+\mu(1-\eta)$.
Since $|\lambda|=|\mu|=1$,
\[
 |h(t)|
 \le\eta(t)+(1-\eta(t))=1
 \qquad(t\in[0,1]).
\]
Moreover, $h(x)=\lambda$ and $h(y)=\mu$.
Thus $\supnorm h=1$, so $h\in S_{\F}$ and $h\in F_{x,\lambda}\cap F_{y,\mu}$.
Therefore the two fibers are not disjoint.

Since any two phase fibers associated with distinct points have nonempty intersection, every family of pairwise disjoint phase fibers containing at least two members must have a common point label.
For each fixed $x$, the fibers $F_{x,\lambda}$ are pairwise disjoint as $\lambda$ varies over $\mathbb T_{\F}$.
Consequently, the maximal families of pairwise disjoint phase fibers are exactly
\[
 \{F_{x,\lambda}:
   \lambda\in\mathbb T_{\F}\},
 \qquad x\in[0,1].
\]
Since $\Delta$ is bijective and permutes the phase fibers, it also permutes these maximal families.
Consequently, there is a bijection $\tau:[0,1]\to[0,1]$ such that $\Delta$ maps the family associated with $x$ onto the family associated with $\tau(x)$.

For each $x$, the restriction of this permutation induces a bijection $\alpha_x:\mathbb T_{\F} \to\mathbb T_{\F}$.
Hence
\[
 \Delta(F_{x,\lambda})=F_{\tau(x),\alpha_x(\lambda)}
 \qquad(x\in[0,1],\ \lambda\in\mathbb T_{\F}),
\]
proving \eqref{eq:boundarypermutation}.

It remains to determine the form of $\alpha_x$.
For $\lambda,\mu\in\mathbb T_{\F}$, evaluation at $x$ gives
\[
 \rp(f,g)\ge|\lambda-\mu|
 \qquad
 (f\in F_{x,\lambda},\
  g\in F_{x,\mu}).
\]
On the other hand, the constant functions $\lambda\one$ and $\mu\one$ belong to the respective fibers and satisfy
\[
 \rp(\lambda\one,\mu\one)
 =|\lambda-\mu|,
\]
because $\qn(a\one)=|a|$ and $\rp(a\one,0)=|a|$ for every $a\in\F$.
Therefore,
\begin{equation}\label{eq:boundarychord}
 \inf\{\rp(f,g):
 f\in F_{x,\lambda},\
 g\in F_{x,\mu}\}
 =|\lambda-\mu|.
\end{equation}

Since $\Delta$ maps each phase fiber onto its corresponding image fiber and preserves distances, \eqref{eq:boundarychord} implies
\[
 |\alpha_x(\lambda)-\alpha_x(\mu)| = |\lambda-\mu|
 \qquad(\lambda,\mu\in\mathbb T_{\F}).
\]
Thus $\alpha_x$ is a surjective Euclidean isometry of $\mathbb T_{\F}$.

If $\F=\R$, then $\mathbb T_{\R}=\{-1,1\}$.
Writing $v_x=\alpha_x(1)$, we obtain $\alpha_x(\lambda)=v_x\lambda$, where $v_x\in\{-1,1\}$.
Suppose now that $\F=\C$.
Put $v_x=\alpha_x(1)$ and define
\[
 \beta_x(\lambda)
 :=\overline{v_x}\alpha_x(\lambda)
 \qquad(\lambda\in\mathbb T_{\C}).
\]
Then $\beta_x$ is a circle isometry fixing $1$.
Distance preservation gives $\beta_x(-1)=-1$ and $\beta_x(i)=i$ or $-i$.
For any $\lambda\in\mathbb T_{\C}$, the equality $|\beta_x(\lambda)-1|=|\lambda-1|$ implies $\operatorname{Re}\beta_x(\lambda) = \operatorname{Re}\lambda$.
If $\beta_x(i)=i$, preservation of distances to $i$ also gives $\operatorname{Im}\beta_x(\lambda)=\operatorname{Im}\lambda$, and hence $\beta_x(\lambda)=\lambda$.
If $\beta_x(i)=-i$, the same argument gives $\operatorname{Im}\beta_x(\lambda)=-\operatorname{Im}\lambda$, so $\beta_x(\lambda)=\overline\lambda$.

Consequently, $\alpha_x(\lambda)=v_x\lambda$ or $v_x\overline\lambda$.
In either scalar field, the asserted formula extends $\alpha_x$ to a $\R$-linear Euclidean isometry $\widetilde\alpha_x:\F\to\F$ fixing zero.
\end{proof}

\begin{lemma}\label{lem:boundaryevaluation}
For $f\in S_{\F}$, $x\in[0,1]$, and $\lambda\in\mathbb T_{\F}$, put
\[
 \mathfrak{f}_{x,\lambda}(f)=\inf_{g\in F_{x,\lambda}}\rp(f,g).
\]
Then
\begin{equation}\label{eq:boundaryinfimum}
 \begin{cases}
 \mathfrak{f}_{x,\lambda}(f)=|\lambda-f(x)|,&|\lambda-f(x)|\ge1/2,\\
 \mathfrak{f}_{x,\lambda}(f)<1/2,&|\lambda-f(x)|<1/2.
 \end{cases}
\end{equation}
With the notation of Lemma~\ref{lem:boundaryfibers},
\begin{equation}\label{eq:fullpoints}
 (\Delta f)(\tau(x))=\widetilde\alpha_x(f(x))
 \qquad(f\in S_{\F},\ x\in[0,1]).
\end{equation}
\end{lemma}

\begin{proof}
We first establish
\[
 \inf_{g\in F_{x,\lambda}}
 \supnorm{f-g}=|\lambda-f(x)|.
\]
Evaluation at $x$ gives the lower bound.
Given $\delta>0$, choose a smooth cutoff $0\le\eta\le1$ with $\eta(x)=1$ supported where $|\lambda-f(t)|<|\lambda-f(x)|+\delta$.
Then $g=(1-\eta)f+\lambda\eta$ takes values in the closed unit ball of $\F$, belongs to $F_{x,\lambda}$, and satisfies
\[
 \supnorm{f-g}\le|\lambda-f(x)|+\delta.
\]
This proves the equality.
Proposition~\ref{prop:large} and $\Hp<1/2$ then give~\eqref{eq:boundaryinfimum}, including the boundary case.

Fix $f$ and $x$, and write
\[
 z=f(x),\qquad
 w=\widetilde\alpha_x^{-1}((\Delta f)(\tau(x))).
\]
Both values lie in the closed unit ball of $\F$.
Surjectivity on fibers gives
\[
 \mathfrak{f}_{\tau(x),\alpha_x(\lambda)}(\Delta f)=\mathfrak{f}_{x,\lambda}(f).
\]
Whenever $|\lambda-z|>1/2$, the common infimum is greater than $1/2$.
Equation~\eqref{eq:boundaryinfimum} in the source and target therefore yields
\begin{equation}\label{eq:boundarydistances}
 |\lambda-w|=|\lambda-z|.
\end{equation}

If $\F=\R$, choose $\lambda\in\{-1,1\}$ with $|\lambda-z|\ge1$.
For fixed such $\lambda$, the map $t\mapsto|\lambda-t|=1-\lambda t$ is injective on $[-1,1]$.
Thus~\eqref{eq:boundarydistances} implies $w=z$.
If $\F=\C$, there is a nonempty open arc of the unit circle on which $|\lambda-z|>1/2$.
For $z\ne0$, take an arc about $-z/|z|$; for $z=0$, any arc works.
Equation~\eqref{eq:boundarydistances} holds on that arc.
Three distinct points of the arc are noncollinear, and equality of distances to those three points forces $w=z$.
Explicitly, subtracting their squared-distance equations gives two independent real linear equations for $w-z$.
This proves~\eqref{eq:fullpoints} over both fields.
\end{proof}

\begin{proposition}\label{prop:fullcomposition}
There exist a smooth function $u:[0,1]\to\mathbb{T}_{\F}$, a smooth diffeomorphism $\sigma:[0,1]\to[0,1]$, and a map $\kappa:\F\to\F$ such that
\[
 \Delta f=u\,\kappa(f\circ\sigma)\qquad(f\in S_{\F}).
\]
The map $\kappa$ is the identity in the real case and either the identity or conjugation in the complex case.
When $\F=\R$, the function $u$ is a constant sign.
The same choice of $\kappa$ applies throughout $[0,1]$ in the complex case.
\end{proposition}
\begin{proof}
Put $u=\Delta\one$.
By~\eqref{eq:fullpoints}, $u(\tau(x))=v_x$, so $u$ is smooth and has modulus one.
For $\F=\R$, connectedness of the interval implies that $u$ is a constant sign, and $\widetilde\alpha_x(z)=v_xz$.
For $\F=\C$, the smooth function
\[
 y\longmapsto\frac{\Delta(i\one)(y)}{i u(y)}
\]
takes values in $\{-1,1\}$ according to the orientation of $\widetilde\alpha_{\tau^{-1}(y)}$.
Connectedness makes it constant, so there is one choice of $\kappa$, either the identity or conjugation, valid at every point.
With $\sigma=\tau^{-1}$, equation~\eqref{eq:fullpoints} gives the claimed formula in either field.

Apply the formula to the real coordinate function $\iota(t)=t\in S_{\F}$.
Since $\kappa$ fixes real numbers,
\[
 \sigma=\frac{\Delta\iota}{\Delta\one}.
\]
Thus $\sigma$ is smooth.
Applying the same reconstruction to $\Delta^{-1}$ yields
\[
 \tau=\frac{\Delta^{-1}\iota}{\Delta^{-1}\one},
\]
so the inverse map is smooth as well.
The chain rule for $\sigma\circ\tau=\id$ shows that their derivatives do not vanish, with the endpoint derivatives understood as continuous extensions.
Hence $\sigma$ is a smooth interval diffeomorphism.
\end{proof}

The representation supplies a $\R$-linear supremum-norm isometry on the ambient space.
The following lemma proves that it also preserves the given metric; no rigidity of the weight or map is needed for this step.

\begin{lemma}\label{lem:ambienttransfer}
Let $X$ be a real normed space of dimension at least two, and let $S_X$ be its unit sphere.
Then
\[
 S_X-S_X=\{h\in X:\norm h\le2\}.
\]
Suppose in addition that $d$ is a translation-invariant metric on $X$ and that, for some $0<c\le2$,
\[
 d(h,0)=\norm h\qquad\hbox{whenever }\norm h\ge c.
\]
If $A:X\to X$ is a $\R$-linear norm isometry and its restriction to $S_X$ preserves $d$, then $A$ preserves $d$ on all of $X$.
\end{lemma}
\begin{proof}
The triangle inequality proves one inclusion in the difference formula.
For $0<\norm h\le2$, set $a=h/\norm h$ and choose $b$ linearly independent of $a$.
The path
\[
 v(t)=\frac{(\cos t)a+(\sin t)b}{\norm{(\cos t)a+(\sin t)b}},\qquad 0\le t\le\pi,
\]
is well defined and joins $a$ to $-a$ in $S_X$.
The continuous function $t\mapsto\norm{v(t)-h}$ starts at $|1-\norm h|\le1$ and ends at $1+\norm h>1$.
By the intermediate value theorem, $t\mapsto\norm{v(t)-h}$ takes the value one.
At such a point, $h=v(t)-(v(t)-h)$ is a difference of two unit vectors.
For $h=0$, take the same unit vector twice.

If $\norm h<c$, write $h=f-g$ with $f,g\in S_X$.
Real linearity and translation invariance give
\[
 d(Ah,0)=d(Af,Ag)=d(f,g)=d(h,0).
\]
If $\norm h\ge c$, norm preservation and the assumed large-distance identity give the same equality.
Applying this to arbitrary differences proves the assertion.
\end{proof}

\begin{proposition}\label{prop:fullambient}
Every surjective $\rp$-isometry $\Delta:S_\F\to S_\F$ has a unique $\R$-linear surjective $\rp$-isometric
extension to the ambient space.
This extension preserves the supremum norm and has the weighted composition form of Proposition~\ref{prop:fullcomposition}.
\end{proposition}
\begin{proof}
Take the representation in Proposition~\ref{prop:fullcomposition} and define
\[
 \widetilde{\Delta}f:=u\,\kappa(f\circ\sigma)\qquad(f\in C^\infty([0,1], \F)).
\]
The map $\widetilde{\Delta}$ is $\R$-linear, preserves the supremum norm, and maps smooth functions to smooth functions.
It is surjective, with inverse
\[
 \widetilde{\Delta}^{-1}g=\kappa\bigl((g/u)\circ\sigma^{-1}\bigr).
\]
Here $\kappa^2=\id$ in either field.
For $\F=\C$, $\widetilde{\Delta}$ is $\C$-linear or conjugate linear according to the choice of $\kappa$.
For $\F=\R$, the weight is a constant sign by Proposition~\ref{prop:fullcomposition}.

Regard $ C^\infty([0,1], \F)$ as a real normed space with the supremum norm.
Its real dimension is at least two, and completeness in this norm is not needed.
Proposition~\ref{prop:large} supplies the large-distance hypothesis of Lemma~\ref{lem:ambienttransfer} with $c=1/2$.
Since $\widetilde{\Delta}|_{S_{\F}}=\Delta$, that lemma proves that $\widetilde{\Delta}$ preserves $\rp$ on all of $ C^\infty([0,1], \F)$.
For any $\R$-linear extension $T$ and any nonzero $f$,
\[
 Tf=\supnorm f\,\Delta(f/\supnorm f)=\widetilde{\Delta}f.
\]
Also $T0=0=\widetilde{\Delta}0$, so the extension is unique and has the asserted radial formula.
\end{proof}

The preceding proposition establishes the extension without requiring the weight to be constant.
We next show that preservation of $\rp$ forces the weight to be constant and the composition map to be admissible.

\begin{lemma}[Rigidity of complex weighted composition isometries]\label{lem:complexrigidity}
Let $u:[0,1]\to\mathbb T_{\C}$ be smooth and $\sigma:[0,1]\to[0,1]$ a smooth diffeomorphism.
The operator $Wf=u(f\circ\sigma)$ is a surjective $\rp$-isometry of $C^\infty([0,1],\C)$ if and only if $u$ is constant and $\sigma$ is admissible.
Equivalently, these conditions characterize when $W|_{S_\C}$ is a surjective $\rp$-isometry.
\end{lemma}
\begin{proof}
The operator $W$ is a surjective $\C$-linear supremum-norm isometry.
Lemma~\ref{lem:ambienttransfer} and Proposition~\ref{prop:large} show that its restriction to $S_\C$ is a metric isometry if and only if $W$ is an ambient metric isometry.
Assume the latter, and put
\[
 \mathcal G=\{h\in C^\infty([0,1],\C):q_{1,p}(h)>6\supnorm h\}.
\]
For $h\in\mathcal G$, set $t=2/q_{1,p}(h)$.
Then $t\supnorm h<1/3$, the first-order contribution to $\rp(th,0)$ is $\tfrac12\phi(2)=1/3$, and every contribution with $n\ge2$ is at most $1/4$.
Thus $\rp(th,0)=1/3$.
Since $W$ preserves both metrics and is linear, the same isolation applies to $tWh$ and gives
\begin{equation}\label{eq:weightedfirst}
 q_{1,p}(u(h\circ\sigma))=q_{1,p}(h)\qquad(h\in\mathcal G).
\end{equation}

\emph{First, $\sigma$ is the identity or reflection.}
Set $h_k(x)=1-\cos(2\pi kx)$ for positive integers $k$.
Then $\supnorm{h_k}=2$ and $\supnorm{h_k'}=2\pi k$, so $h_k\in\mathcal G$ for large $k$.
With $M_k$ as in Lemma~\ref{lem:oscillation}, differentiation gives
\[
 (u(h_k\circ\sigma))'=u'(h_k\circ\sigma)
       +2\pi k u\sin(2\pi k\sigma)\sigma',
\]
and hence
\[
 \left|\frac{\supnorm{(u(h_k\circ\sigma))'}}{2\pi k}-M_k\right|
 \le\frac{\supnorm{u'}}{\pi k}.
\]
The anchored zeroth-order values have modulus at most two.
For either finite $p$ or $p=\infty$, equation~\eqref{eq:weightedfirst}, divided by $2\pi k$, therefore yields $M_k\to1$.
Lemma~\ref{lem:oscillation} implies $\supnorm{\sigma'}=1$.
As $\sigma$ is a monotone bijection,
\[
 \int_0^1|\sigma'(x)|\,dx=|\sigma(1)-\sigma(0)|=1.
\]
Continuity forces $|\sigma'|\equiv1$, and consequently $\sigma=\id$ or $\sigma=r$.

\emph{Second, $u$ is constant.}
Since $h_k\circ r=h_k$, equation~\eqref{eq:weightedfirst} now gives
\begin{equation}\label{eq:weighttest}
 q_{1,p}(uh_k)=q_{1,p}(h_k)\qquad\text{for all sufficiently large }k.
\end{equation}
Differentiating $|u|^2=1$ gives $\operatorname{Re}(\overline u u')=0$.
For every $\R$-valued smooth $h$, it follows that
\[
 |(uh)'(x)|^2=|h'(x)|^2+|u'(x)|^2h(x)^2.
\]
If $u'\not\equiv0$, there are a nondegenerate closed interval $I\subset(0,1)$ and $\delta>0$ with $|u'|\ge\delta$ on $I$.
For all sufficiently large $k$, choose $x_k\in I$ with $|\sin(2\pi kx_k)|=1$.
Then $h_k(x_k)=1$ and $|h_k'(x_k)|=2\pi k$, so
\[
 \supnorm{(uh_k)'}^2\ge(2\pi k)^2+\delta^2>(2\pi k)^2.
\]
The anchored values of $uh_k$ and $h_k$ have the same modulus.
For finite $p$, the strict derivative increase contradicts~\eqref{eq:weighttest} directly.
For $p=\infty$, use $|h_k(x_0)|\le2<2\pi k$, so that $q_{1,\infty}(h_k)=2\pi k$; the same contradiction follows.
Thus $u'=0$, and $u=\lambda$ for a constant $|\lambda|=1$.

\emph{Third, reflection must be admissible.}
Suppose $p<\infty$, $\sigma=r$, and $x_0\ne1/2$.
Choose a nonzero smooth function $\chi\ge0$ supported in an open subinterval of $(0,1)$ disjoint from $\{x_0,1-x_0\}$, and set
\[
 g_k(x)=(x-x_0)^2+\chi(x)(1+\sin(2\pi kx)).
\]
Their supremum norms are uniformly bounded.
On an interval where $\chi\ge\eta>0$, points with $|\cos(2\pi kx)|=1$ give
\[
 \supnorm{g_k'}\ge2\pi k\eta-2-2\supnorm{\chi'}\longrightarrow\infty.
\]
Thus $g_k\in\mathcal G$ for large $k$, whereas
$g_k(x_0)=0$ and $g_k(1-x_0)=(1-2x_0)^2>0$.
Equation~\eqref{eq:weightedfirst} would imply
\[
 \supnorm{g_k'}^p=q_{1,p}(g_k)^p
 =q_{1,p}(\lambda(g_k\circ r))^p
 =|g_k(1-x_0)|^p+\supnorm{g_k'}^p,
\]
a contradiction.
Hence reflection for finite $p$ requires $x_0=1/2$.

Conversely, multiplication by a constant unimodular scalar preserves all the defining norms.
Composition by an admissible map preserves $\rp$ by Propositions~\ref{prop:finite} and~\ref{prop:anchorfree}.
These operators are surjective, which proves sufficiency.
\end{proof}

\begin{proof}[\textbf{Proof of Theorem~\ref{thm:full}}]
Proposition~\ref{prop:fullambient} supplies the unique $\R$-linear isometric extension $\widetilde{\Delta}f=u\,\kappa(f\circ\sigma)$.
If $\F=\C$ and $\kappa$ is the identity, Lemma~\ref{lem:complexrigidity} gives a constant unimodular weight and an admissible map.
If $\kappa$ is conjugation, compose $\widetilde{\Delta}$ on the right with $Jf=\overline f$.
Conjugation preserves all the defining norms, and $(\widetilde{\Delta}\circ J)f=u(f\circ\sigma)$, so Lemma~\ref{lem:complexrigidity} applies.
For $\F=\R$, the weight is a constant sign $\varepsilon$.
The operator $\varepsilon \widetilde{\Delta}$ maps $\Sp$ onto itself, so Theorem~\ref{thm:main} gives an admissible map $\sigma_0$ there.
Evaluation of both composition formulas at $\iota(t)=t$ yields $\sigma=\sigma_0$.

Conversely, multiplication by a constant unimodular scalar, composition by an admissible map, and complex conjugation are all surjective $\rp$-isometries of the ambient space that preserve $S_\F$.
Their compositions therefore induce surjective $\rp$-isometries of $S_\F$.
\end{proof}

\begin{corollary}\label{cor:fullconsequences}
For $\F\in\{\R,\C\}$, every surjective $\rp$-isometry $\Delta:S_{\F}\to S_{\F}$ satisfies
\[
 \Delta(-f)=-\Delta f,\qquad
 \supnorm{\Delta f-\Delta g}=\supnorm{f-g}\qquad(f,g\in S_{\F}).
\]
Its extension is the radial map
\begin{equation}\label{eq:fullradial}
 \widetilde{\Delta}f=\begin{cases}
 \supnorm f\,\Delta(f/\supnorm f),&f\ne0,\\
 0,&f=0.
 \end{cases}
\end{equation}
Moreover, $\widetilde{\Delta}$ is the unique surjective global $\rp$-isometry extending $\Delta$, even when linearity of a competing global extension is not assumed.
\end{corollary}
\begin{proof}
The $\R$-linear, supremum-norm-preserving representations in Theorem~\ref{thm:full} prove antipodal and supremum-distance preservation.
Real linearity gives~\eqref{eq:fullradial}.
If $T: C^\infty([0,1], \F)\to C^\infty([0,1], \F)$ is another surjective global isometric extension, Corollary~\ref{cor:affine} gives $T(f)=b+L(f)$ with $L$ real linear.
For $f\in S_{\F}$, antipodal preservation yields
\[
 2b=T(f)+T(-f)=\Delta f+\Delta(-f)=0.
\]
Thus $T$ is $\R$-linear, and Theorem~\ref{thm:full} gives $T=\widetilde{\Delta}$.
\end{proof}

\begin{remark}
Formula~\eqref{eq:fullradial} determines the only possible $\R$-linear extension, but its additivity follows here from the composition representations, and ambient metric preservation over both value fields follows from Lemma~\ref{lem:ambienttransfer}.
They are not consequences of radial normalization alone, since $\rp$ is not a homogeneous metric.
\end{remark}
The full-sphere classification, together with global affinity, also determines all surjective isometries of the ambient spaces.

\begin{corollary}\label{cor:global}
Let $\F\in\{\R,\C\}$.
A map $T:C^\infty([0,1],\F)\to C^\infty([0,1],\F)$ is a surjective $\rp$-isometry if and only if
\begin{equation}\label{eq:global}
 T(f)=b+\lambda\,\kappa(f\circ\sigma)
 \qquad(f\in C^\infty([0,1],\F)),
\end{equation}
where $b\in C^\infty([0,1],\F)$, $\lambda\in\mathbb T_{\F}$ is constant, $\sigma$ is admissible, and $\kappa$ is the identity when $\F=\R$ and either the identity or complex conjugation when $\F=\C$.
\end{corollary}
\begin{proof}
Suppose that $T$ is a surjective $\rp$-isometry.
Put $b=T(0)$ and $L(f)=T(f)-b$.
By translation invariance, $L$ is a surjective $\rp$-isometry fixing zero.
Theorem~\ref{thm:recovery} and Corollary~\ref{cor:affine} show that $L$ is $\R$-linear and preserves the supremum norm.
Consequently, its restriction to $S_{\F}$ is a surjective $\rp$-isometry.
Theorem~\ref{thm:full}, together with uniqueness of the $\R$-linear extension, gives
\[
 L(f)=\lambda\,\kappa(f\circ\sigma)
 \qquad(f\in C^\infty([0,1],\F))
\]
with the stated parameters.

Conversely, Theorem~\ref{thm:full} shows that $f\mapsto\lambda\,\kappa(f\circ\sigma)$ is a surjective $\rp$-isometry for every permitted choice of parameters.
Since translations preserve $\rp$, the map in \eqref{eq:global} is also a surjective $\rp$-isometry.
\end{proof}

\section*{Acknowledgments}
The author acknowledges the use of OpenAI's ChatGPT as an AI-assisted tool in preparing this manuscript, including assistance with revising mathematical arguments, improving the exposition, and editing the \LaTeX{} source.
The author independently checked the mathematical arguments and references, revised the generated material, and takes full responsibility for the final manuscript.

\appendix
\section{Metric ball formulas and Minkowski functionals}\label{app:balls}
We write $B_d(a,r)$ and $\overline B_d(a,r)$ for the open and closed metric balls with respect to the induced metric $d$.
\begin{proposition}[Metric balls for the finite-order metrics]\label{prop:finite-order-balls}
Let $m\ge0$, $1\le p\le\infty$, and $r>0$, and put $I(r):=\{n\ge1:2^{-n}>r\}$.
Then $I(r)$ is finite, and the closed and open balls of $\rho_{m,p}$ centered at zero are
\begin{align*}
 \overline B_{\rho_{m,p}}(0,r)
 &=\{f:q_{m,p}(f)\le r\}\cap
   \bigcap_{n\in I(r)}
   \left\{f:\qn(f)\le
   \frac{r}{2^{-n}-r}\right\},\\
 B_{\rho_{m,p}}(0,r)
 &=\{f:q_{m,p}(f)<r\}\cap
   \bigcap_{n\in I(r)}
   \left\{f:\qn(f)<
   \frac{r}{2^{-n}-r}\right\}.
\end{align*}
In particular, every $\rho_{m,p}$-ball is convex, and every positive-radius ball centered at zero is balanced and absorbing.
For $r\ge1/2$, the balls coincide with the corresponding $q_{m,p}$-norm balls.
\end{proposition}
\begin{proof}
Recall that
\[
 \rho_{m,p}(f,0)
 =\max\left\{
 q_{m,p}(f),
 \sup_{n\ge1}2^{-n}
 \frac{\qn(f)}{1+\qn(f)}
 \right\}.
\]
For every $n\ge1$, we have
\[
 0\le
 2^{-n}\frac{\qn(f)}{1+\qn(f)}
 <2^{-n}.
\]
If $2^{-n}\le r$, both the strict and non-strict inequalities against $r$ hold automatically.
If $2^{-n}>r$, then
\[
 2^{-n}\frac{\qn(f)}{1+\qn(f)}
 \le r
 \quad\Longleftrightarrow\quad
 \qn(f)\le\frac{r}{2^{-n}-r},
\]
with the same equivalence for strict inequalities.
This proves the closed-ball formula.

For the open ball, it remains to justify that the pointwise strict inequalities make the supremum strictly less than $r$.
Choose $N$ such that $2^{-(N+1)}<r$.
The finitely many terms with $1\le n\le N$ are each strictly less than $r$, while
\[
 \sup_{n>N}2^{-n}
 \frac{\qn(f)}{1+\qn(f)}
 \le2^{-(N+1)}<r.
\]
Hence their full supremum is strictly less than $r$, proving the open-ball formula.

Each displayed ball is a finite intersection of norm balls, so it is convex, balanced, and absorbing.
Translation invariance then implies convexity of balls about arbitrary centers.
Finally, if $r\ge1/2$, then $I(r)=\varnothing$; the formulas reduce to the open and closed $q_{m,p}$-norm balls.
\end{proof}

\begin{corollary}\label{cor:gauges}
The Minkowski functional of $\overline B_{\rp}(0,r)$ is the continuous norm
\[
 \gamma_r(f)=\max\left\{\frac{q_0(f)}r,\,
          \max_{n\in I(r)}\frac{2^{-n}-r}{r}\qn(f)\right\},
\]
with the inner maximum omitted if $I(r)=\varnothing$.
The family $(\gamma_r)_{r>0}$ generates the smooth topology, but no single $\gamma_r$ does.
For $r\ge1/2$, $\gamma_r=q_0/r$.
\end{corollary}
\begin{proof}
The ball formulas give $\overline B_{\rp}(0,r)=\{\gamma_r\le1\}$ and $B_{\rp}(0,r)=\{\gamma_r<1\}$.
Homogeneity identifies $\gamma_r$ with the Minkowski functional.
For fixed $r$, put $N=\max I(r)$ if this set is nonempty, and $N=0$ otherwise.
Proposition~\ref{prop:norms} shows that $\gamma_r$ is equivalent to $P_N$.

For every $n\ge1$, choose $0<r<2^{-n}$.
Then
\[
 \gamma_r(f)\ge
 \frac{2^{-n}-r}{r}\qn(f).
\]
Conversely, each $\gamma_r$ is a finite maximum of continuous norms $\qn$, together with $q_0/r$.
Hence the family $(\gamma_r)_{r>0}$ generates the smooth topology $\mathcal T$.

No finite-order norm $P_N$ induces the full smooth topology: for $h_k(x)=(2\pi k)^{-(N+1)}\sin(2\pi kx)$, we have $P_N(h_k)\to0$ but $P_{N+1}(h_k)=1$.
Consequently, no single $\gamma_r$ generates $\mathcal T$.
\end{proof}

\begin{remark}
Convexity of the balls does not make $\rp$ a norm metric.
Its small balls impose derivative conditions, whereas its balls of radius at least $1/2$ are supremum norm balls.
They are not all scalar dilates of one neighborhood.
\end{remark}

\section{Anchored-coordinate representation}\label{app:anchored}
\begin{proposition}\label{prop:anchored-coordinates}
Let $\F\in\{\R,\C\}$, $n\ge1$, and $1\le p\le\infty$.
The map
\[
 J_n f=\bigl(f(x_0),f'(x_0),\ldots,f^{(n-1)}(x_0),f^{(n)}\bigr)
\]
is a $\F$-linear surjective isometry from $(C^n([0,1],\F),\qn)$ onto $\ell_p^n(\F)\oplus_p (C([0,1],\F),\norm{\cdot}_\infty)$, where the direct sum carries the maximum norm when $p=\infty$.
Its inverse is given by
\[
 J_n^{-1}(a_0,\ldots,a_{n-1},g)(x)
 =
 \sum_{j=0}^{n-1} a_j\frac{(x-x_0)^j}{j!} 
 + \frac{1}{(n-1)!}\int_{x_0}^x (x-t)^{n-1}g(t)\,dt,\quad x\in[0,1].
\]
Moreover, the restriction of $J_n$ to $(C^\infty([0,1],\F),\qn)$ is a $\F$-linear surjective isometry onto $\ell_p^n(\F)\oplus_p (C^\infty([0,1],\F), \norm{\cdot}_\infty)$.
\end{proposition}

\begin{proof}
The linearity of $J_n$ is immediate.
By the definition of $\qn$,
\[
 \|J_nf\|_{\ell_p^n\oplus_p C}
 =\qn(f)
 \qquad(f\in C^n([0,1],\F)),
\]
so $J_n$ is an isometry.
To prove surjectivity, let $(a_0,\ldots,a_{n-1},g)
\in \ell_p^n(\F)\oplus_p C([0,1],\F)$.
Define
\[
f(x)=
\sum_{j=0}^{n-1}
a_j\frac{(x-x_0)^j}{j!}
+\frac{1}{(n-1)!}
\int_{x_0}^x(x-t)^{n-1}g(t)\,dt\qquad(x\in[0,1]).
\]
Since $g$ is continuous, successive differentiation under the integral sign gives, for $0\le k<n$,
\begin{equation}\label{eq:anchored-inverse-derivatives}
f^{(k)}(x)=
\sum_{j=k}^{n-1}
a_j\frac{(x-x_0)^{j-k}}{(j-k)!}
+\frac{1}{(n-k-1)!}
\int_{x_0}^x
(x-t)^{n-k-1}g(t)\,dt.
\end{equation}
In particular, evaluating at $x=x_0$ yields $f^{(k)}(x_0)=a_k$ for each $0\le k<n$.
For $k=n-1$, the derivative formula~\eqref{eq:anchored-inverse-derivatives} becomes
\[
f^{(n-1)}(x)=a_{n-1}+\int_{x_0}^x g(t)\,dt.
\]
The fundamental theorem of calculus therefore gives $f^{(n)}=g$ on $[0,1]$.
Thus $f\in C^n([0,1],\F)$ and $J_nf=(a_0,\ldots,a_{n-1},g)$, proving surjectivity.

Since $J_n$ is a linear isometry, it is injective, and hence a linear isometric bijection.
The construction above gives its inverse.
Finally, if $g\in C^\infty([0,1],\F)$, then $f^{(n)}=g$ implies that $f\in C^\infty([0,1],\F)$.
Conversely, the $n$th derivative of every smooth function is smooth.
Therefore the restriction of $J_n$ has the stated range and remains a linear isometric bijection.
\end{proof}

\end{document}